\documentclass{commatDV}

\usepackage{graphicx}

%
%
\newcommand {\supplus}{\mathop{{\supset}\llap{\raise
0.5pt\text{\normalfont\small+}\hskip 0.5pt}}}

\newcommand {\subplus}{\mathop{{\subset}\llap{\raise
0.5pt\text{\normalfont\small+}\hskip 0.5pt}}}

\newcommand {\Cee} {{\mathbb C}}

\newcommand {\Kee} {{\mathbb K}}

\newcommand {\Ree} {{\mathbb R}}

\newcommand {\Zee} {{\mathbb Z}}

\newcommand {\fg} {{\mathfrak{g}}} %
\newcommand {\fgl} {{\mathfrak{gl}}} %

\newcommand {\fsl} {{\mathfrak{sl}}}

\newcommand {\fsvect} {{\mathfrak{svect}}}

\newcommand {\fvect} {{\mathfrak{vect}}} %

\renewcommand {\cal} {\mathcal}

\newcommand {\cD} {{\cal D}}

\newcommand {\cF} {{\cal F}}

\newcommand {\cI} {{\cal I}}

\newcommand {\cL} {{\cal L}}

\newcommand {\cU} {{\cal U}}

%
%

\def \opname#1#2%
 {\expandafter\newcommand \csname #1\endcsname {{\mathop{#2}\nolimits}}}


\newcommand{\rmname}[1]
 {\expandafter\newcommand \csname #1\endcsname {{\operatorname{#1}}}}

\newcommand{\rmnameii}[2]
 {\expandafter\newcommand \csname #1\endcsname {{\operatorname{#2}}}}

\rmname{act}
\rmname{Ad}
\rmname{Add}
\rmname{ad}
\rmname{Alt}
\rmname{alt}
\rmname{Ann}
\rmname{antidiag}
\rmname{Ber}
\rmname{ber}
\rmname{Br}
\rmname{card}
\rmname{ch}
\rmname{Char}
\rmname{cem}
\rmname{cj}
\rmname{Cliff}
\rmname{cntr}
\rmname{codim}
\rmname{coind}
\rmname{const}
\rmname{col}
\rmname{cork}
\rmname{cpr}
\rmname{diag}
\rmnameii{Div}{div}
\rmname{Def}
\rmname{Der}
\rmname{Dim}
\rmname{End}
\rmname{Even}
\rmname{Ext}
\rmname{gr}
\rmname{Hom}
\rmname{HT}
\rmnameii{Ht}{ht}
\rmname{hwt}
\rmname{Id}
\rmname{id}
\rmname{ind}
\rmname{Ind}
\rmname{Inf}
\rmname{irr}
\rmname{Le}
\rmname{Lie}
\rmname{lwt}
\rmname{mult}
\rmname{Mor}
\rmname{nm}
\rmname{Ob}
\rmname{Odd}
\rmname{Osc}
\rmname{per}
\rmname{Pic}
\rmname{pr}
\rmname{pro}
\rmname{Prime}
\rmname{Proj}
\rmname{prt}
\rmname{pt}
\rmname{Q}
\rmname{qet}
\rmname{qtr}
\rmname{rd}
\rmname{rk}
\rmname{row}
\rmname{Res}
\rmname{salt}
\rmname{Sch}
\rmname{SBr}
\rmname{scalar}
\rmname{Ser}
\rmname{sign}
\rmname{Smbl}
\rmname{spin}
\rmname{ssym}
\rmname{str}
\rmname{st}
\rmname{sgn}
\rmname{sq}
\rmname{symm}
\rmname{supp}
\rmname{Supp}
\rmname{St}
\rmname{Spec}
\rmname{Specm}
\rmname{Spm}
\rmname{tr}
\rmname{vpt}
\rmname{weyl}
\rmname{Weyl}
\rmname{Witt}

\opname{vvol} {{v\hspace{-0.1ex}o\hspace{-0.02ex}l\/}}
\opname{pnt} {\text{\normalfont pt}}
\opname{Span} {{Span}}
\opname{slim} {\overline{\lim}}
\opname{Vol} {{V\hspace{-0.55ex}o\hspace{-0.02ex}l\/}}
\opname{QVol} {{Q\hspace{-0.3ex}V\hspace{-0.55ex}o\hspace{-0.02ex}l\/}}
\opname{PoVol}{{P\hspace{-0.35ex}o\hspace{-0.25ex}V\hspace{-0.55ex}o\hspace{-0.02ex}l\/}}
\opname{BVol} {{B\hspace{-0.2ex}V\hspace{-0.55ex}o\hspace{-0.02ex}l\/}}
\opname{Par} {{P\hspace{-0.3ex}a\hspace{-0.05ex}r\/}}

%
%

\rmname{Mat}
\rmname{Bil}
\rmname{Diff}
\rmname{Ker}
\rmname{Herm}
\rmname{Coker}
\rmname{Conn}
\rmname{Covect}
\rmname{Vect}
\rmname{Int}

\rmnameii {IM} {Im}
\rmnameii {RE} {Re}

\opname{Aut} {{A\hspace{-0.2ex}u\hspace{-0.1ex}t\/}}
\opname{GL} {{G\hspace{-0.3ex}L}}
\opname{SL} {{S\hspace{-0.3ex}L}}
\opname{Exp} {{E\hspace{-0.2ex}x\hspace{-0.1ex}p\/}}
\opname{GQ} {{G\hspace{-0.2ex}Q}}
\opname{OSp} {{O\hspace{-0.25ex}S\hspace{-0.15ex}p\/}}
\opname{Out} {{O\hspace{-0.25ex}u\hspace{-0.15ex}t\/}}
\opname{Spp} {{S\hspace{-0.2ex}p\/}}
\opname{SpO} {{S\hspace{-0.2ex}p\hspace{-0.02ex}O\/}}
\opname{Pe} {{P\hspace{-0.25ex}e\/}}
\opname{SPe} {{S\hspace{-0.25ex}P\hspace{-0.25ex}e\/}}
\opname{Spin} {{S\hspace{-0.25ex}p\hspace{-0.05ex}i\hspace{-0.1ex}n\/}}
\opname{Iso} {{I\hspace{-0.25ex}s\hspace{-0.1ex}o\/}}
\opname{SSPe} {{S\hspace{-0.25ex}S\hspace{-0.15ex}P\hspace{-0.25ex}e\/}}
\opname{PeU} {{P\hspace{-0.25ex}e\hspace{-0.1ex}U\/}}
\opname{QU} {{Q\hspace{-0.15ex}U\/}}
\opname{U} {{U\/}}

\opname{cGQ} {{\cal G \hspace{-0.2em} Q \/}}
\opname{cSL} {{\cal S \hspace{-0.2em} L \/}}
\opname{cGL} {{\cal G \hspace{-0.2em} L \/}}
\opname{cOSp} {{\cal O \hspace{-0.2em} S \hspace{-0.3em} \it p\/}}
\opname{cPe} {{\cal P \hspace{-1.5pt} \it e\/}}
\opname{cVect} {{\cal V \hspace{-1.5pt} \it e\hspace{-0.1ex}c\hspace{-0.1ex}t\/}}
\opname{cVol} {{\cal V \hspace{-1.5pt} \it o\hspace{-0.1ex}l\/}}
\opname{cAut} {{\cal A \hspace{-0.2em} \it u\hspace{-0.1em}t\/}}
\opname{cCovect} {{\cal C \hspace{-1.5pt}
 \it o\hspace{-0.1ex}v\hspace{-0.1ex}e\hspace{-0.1ex}c\hspace{-0.1ex}t\/}}
\opname{CW} {{C\hspace{-0.15ex}W}}

%
%

\opname {Ab} {{\sf Ab}}
\opname {Alg} {{\sf Alg}}
\opname {ASch} {{\sf Aff\;Sch}}
\opname {Funct} {{\sf Funct}}
\opname {Gr} {{\sf Gr}}
\opname {Grf} {{{\sf Gr}_f}}
\opname {Mods} {{\sf mods}}
\opname {Rings} {{\sf Rings}}
\opname {Salg} {{\sf Salg}}
\opname {Sets} {{\sf Sets}}
\opname {SSMan} {{\sf SMan}}
\opname {Top} {{\sf Top}}
\opname {Vebun} {{\sf Vebun}}

%
%

\newcommand {\tto} {\longrightarrow}

\newcommand {\pderf}[2] {{\frac{\partial {#1}}{\partial {#2}}}}

%
%

\newcommand {\bcdot} {{\mathbin{\text{\raise.4ex\hbox{\bf.}}}}} 

%
%

%
%

\newcommand {\ssecfont} {\normalfont\bf}

\newcommand {\ssec}{\subsection*}

\newcommand {\ssbegin}[1]
  {\refstepcounter{subsection}
  \ea\def\csname the#1\endcsname {{\ssecfont \thesubsection}%
  }%
   \begin{#1}}

\newcommand {\sssbegin}[1]
  {\refstepcounter{subsubsection}
  \ea\def \csname the#1\endcsname {{\ssecfont \thesubsubsection}%
  }%
   \begin{#1}}

   \newcommand {\parbegin}[1]
  {\refstepcounter{paragraph}
  \ea\def \csname the#1\endcsname {{\ssecfont \theparagraph}%
  }%
   \begin{#1}}

\newcommand {\SSbegin}[2]
  {\refstepcounter{subsection}
  \ea\def\csname the#2\endcsname {{\ssecfont #1}%
  }%
   \begin{#2}}

\newcommand {\SSSbegin}[2]
  {\refstepcounter{subsubsection}
  \ea\def \csname the#2\endcsname {{\ssecfont #1}%
  }%
   \begin{#2}}

   \newcommand {\Parbegin}[2]
  {\refstepcounter{paragraph}
  \ea\def \csname the#2\endcsname {{\ssecfont #1}%
  }%
   \begin{#2}}

\setcounter{tocdepth}{2}

\title{Invariant bilinear differential operators}

\author[Pavel Grozman]{\fbox{Pavel Grozman}}

\affiliation{N/A}

\keywords{Representations of Lie algebras, invariant operators.}

\msc{53A55, 17B15 (Primary) 58A99 (Secondary)}

\abstract{Let $M$ be an $n$-dimensional manifold, let $V$ be the space
of a~representation ${\rho: \GL(n)\tto \GL(V)}$. Locally, let $T(V)$ be
the space of sections of the tensor bundle with fiber $V$ over a
sufficiently small open set $U\subset M$, in other words, $T(V)$ is
the space of tensor fields of type $V$ on $U$. In $T(V)$, the group
$\Diff (U)$ of diffeomorphisms of $U$ naturally acts by means of
$\rho$ applied to the Jacobi matrix of the diffeomorphism at the
point.

Here, I give the details of the classification of the $\Diff
(M)$-invariant differential operators $D:T(V_{1})\otimes
T(V_{2})\tto T(V_{3})$ for irreducible fibers with lowest weight.
Up to dualization and the permutation of arguments $T(V_{1})\otimes T(V_{2})\simeq T(V_{2})\otimes T(V_{1})$ a.k.a. ``twist", these operators split into 9 types of operators of order 1, four types of order 2 and 3 types of order 3. The operators of orders 2 and 3 are compositions of 1st order operators with one exception: an indecomposable 3rd order operator which exists only for $n=1$. There are no operators of order $>3$.

Amazingly, almost each 1st order operator determines a~Lie
superalgebra structure on its domain. Moreover, this Lie superalgebra
is almost simple (is a~central extension of a~simple one or contains a
simple ideal of codimention 1).
}

\VOLUME{30}
\NUMBER{3}
\firstpage{129}
\DOI{https://doi.org/10.46298/cm.9341}

\begin{paper}

\tableofcontents

\ssec{Preface}

\addcontentsline{toc}{section}{Preface}

Several years ago Leites told me that the book \cite{KMS}
devoted to natural differential operators contains a~complaint
that the details of the proof of my classification were never
published nor preprinted. Actually, they were deposited to VINITI
and contain not only the proof in the general case but also the
proof in the divergence-free case. However, it was not easy to
retrieve anything from VINITI depositions even during Soviet
period, now they are, it seems, totally inaccessible.

Here are the details of the proof in the general case; the
divergence-free case adds only several compositions of operators, so I skipped
it; for formulation, see \cite{G2}. I also skipped verification
of invariance of several operators, the fact being known, more or
less, to Niujenhuis, and by now can be considered \lq\lq
well-known".

In view of plans greater than the proof of the claims in my Ph.D. thesis' results (these plans are described in \cite{GLS}), Leites
arranged translation and typing of the proof. 
Here is this version. 

I want to warn the reader. There are two detailed expositions of the
proof (a draft and the final version of my thesis). In the final
version, the proof is absolutely correct, but it is written very
succinctly and proof modified (as compared with the draft versions)
with the peculiarities of the case under study being taken into
account. As a~result, the proof is shorter (which was my goal), but is difficult to
generalize to other algebras or to infinite-dimensional fibers.

Stockholm, 2006.

\section*{Introduction}

\addcontentsline{toc}{section}{Introduction}

\ssec{1. Invariant operators: an overview} By {\it invariant
operators} we will mean operators acting in the spaces of tensor
fields (or sections of other types of vector bundles) which have the
same form in any (curvilinear) coordinate system on the fixed manifold
$M$.

The importance of such operators became manifest after discovery of
the relativity theory. Indeed, according to equivalence principle,
the motion of a~body in the gravitational field is equivalent to the
motion in the absence of the field but in a~non-inertial coordinate
system, with curvilinear coordinates if the gravitational field is
non-homogeneous. Thanks to Einstein equations, the action of the
gravitational field on bodies is expressed via the metric of the
space. Invariance of the Einstein equations is a~mathematical
formulation of the equivalence principle.

Similarly, invariant operators should always appear whenever there
exists either a~relation between tensor fields (or sections of vector
bundles depending on higher jets of the diffeomorphism group), or a
condition on a~tensor field, or an algebraic structure, etc., that do
not vary under the changes of coordinates.

Examples: Lie algebra structure on the space of vector fields, the
Stokes formula, the equation of a~geodesic curve, condition for a
local rectifiability of a~pair of vector fields, condition for local
integrability of a~distribution, etc., or, if we confine ourselves to
analytic coordinates only, Cauchy--Riemann equations, etc.

By \textit{tensor fields} we will mean sections of the bundle
\[
E^p_{q}(M)=(TM)^{\otimes p}\otimes (TM^*)^{\otimes q}. 
\]
By
\textit{$\lambda$-densities} we mean sections of $\Vol^{\lambda}=(\Omega^{\dim
M})^{\otimes \lambda}$; the space is well-defined as a module over the group of
diffeomorphisms for non-negative
(and by dualizing for all) integer values of $\lambda$, but
infinitesimally, on the level of the ``Lie algebra of the group of
diffeomorphisms'', the action of the Lie algebra  can be defined for any $\lambda$. As we
will show, for such Lie algebra one can take the Lie algebra of
vector fields with polynomial or formal coefficients and the action of
the vector field in $\Vol^{\lambda}$ is just the multiplication by divergence with
factor $\lambda$.

In this paper we consider the unary linear differential operators
$$
\Gamma(M, E^p_{q}\otimes \Vol^{\lambda}(M))\tto
\Gamma(M, E^r_{s}\otimes \Vol^{\mu}(M))
$$
and bilinear differential operators
$$
\Gamma(M, E^{p_{1}}_{q_{1}}\otimes \Vol^{\lambda_{1}}(M))\times
\Gamma(M, E^{p_{2}}_{q_{2}}\otimes \Vol^{\lambda_{2}}(M))\tto
\Gamma(M, E^r_{s}\otimes \Vol^{\mu}(M)).
$$

The simplest linear invariant operator is the differential of a
function
$$
f\mapsto df=\mathop{\sum}\limits_{i=1}^n dx_{i}\pderf{f}{x_{i}}.
$$
The invariance of this operator is one of the fundamental theorems of
Calculus.

A generalization of this operator is the exterior differential of
differential forms
$$
d: \Omega^{p}\tto \Omega^{p+1}.
$$
It turns out that
$$
d\quad \text{{\sl is the only linear invariant differential operator of nonzero order}}
$$
acting in the spaces of tensor fields with irreducible fibers. This
was proven for more and more general tensors: for differential forms
(\cite{P}, 1959), for covariant tensor fields (\cite{L}, 1973) and,
finally, for general tensors independently and by different methods by
Rudakov \cite{R1} in 1973, Terng (Ph.D.
Thesis, 1976, see \cite{T}), and Kirillov \cite{Ki1} in 1977.

Consider bilinear operators. Historically, the first and most known
first order differential operator is the {\it Lie derivative}
$$
L: \Gamma(M, TM)\times
\Gamma(M, E^{p}_{q}\otimes \Vol^{\lambda}(M))\tto
\Gamma(M, E^{p}_{q}\otimes \Vol^{\lambda}(M)).
$$
Particular cases of this operator: the bracket of vector fields and
operators representable as compositions of $d$ and zero order
operators.

In the first half of XX century, after works by Einstein and Hilbert on
general relativity, researchers started a~systematic search of
invariant operators. Veblen explicitly formulated the problem at the
1928 Mathematical Congress in Bolognia \cite{V}. In 1940 and 1954,
Schouten found two new invariant operators:
$$
\Gamma(M, \Lambda ^kTM)\times
\Gamma(M, \Lambda ^lTM)\tto
\Gamma(M, \Lambda ^{k+l-1}TM)
$$
and
$$
L: \Gamma(M, E^{p}_{q}(M))\times
\Gamma(M, E^{q}_{p}\otimes \Vol^{\lambda}(M))\tto
\Gamma(M, T^*M\otimes \Vol(M)).
$$
called the {\it anti-symmetric} and {\it Lagrangian} concomitant,
respectively.

Schouten also observed that the Poisson bracket can be interpreted,
if one restricts to functions homogenous on fibers, as a~first order invariant
operator ({\it symmetric} concomitant)
$$
P: \Gamma(M, S^kTM)\times
\Gamma(M, S^lTM)\tto
\Gamma(M, S^{k+l-1}TM).
$$
In 1955, a~student of Schouten, Nijenhuis, found one more invariant
operator (the {\it Nijenhuis bracket}) on the space of vector-valued
forms (see \cite{N1})
$$
N: \Gamma(M, TM\otimes\Lambda ^kT^*M)\times
\Gamma(M, TM\otimes\Lambda ^lT^*M)\tto
\Gamma(M, TM\otimes\Lambda ^{k+l}T^*M).
$$
During the next 20 years various applications of these operators were
studied (\cite{Bu}, \cite{FF1}, \cite{FF2}, \cite{N2}, \cite{Tu}).

In 1977--78 in my BS and MS theses, I have completely classified bilinear invariant differential
operators for $\dim M\leq 2$, see \cite{G1}. Three new operators were
found, denoted in what follows $F$, $G$, and $P^*$.

A. Kirillov noticed (see \cite{Ki1}) that by means of the invariant pairing (index $c$
indicates that we consider fields with compact support)
$$
B: \Gamma_c(M, E^{p}_{q}(M))\times \Gamma(M, E^{q}_{p}\otimes
\Vol(M))\tto \Ree
$$
one can define the duals (with respect to the first or second
argument) briefly referred in what follows as the {\it first and second
duals} or 1-dual, $B^{1*}$, and 2-dual, $B^{2*}$, of $B$.

Clearly, if $B$ is invariant, so are its duals. It turned out that
the lagrangian concomitant is dual to the Lie derivative, whereas
the operators dual to the other two Schouten's concomitants and
to the Nijenhuis bracket turned out to be new.

In the same paper, Kirillov generalized to dimensions $\dim M> 2$ the
above mentioned operator $F$:
$$
F: \Gamma(M, \Lambda^{p}T^*M\otimes \Vol^k)\times
\Gamma(M, \Lambda^{q}T^*M\otimes \Vol^l)\tto
\Gamma(M, \Lambda^{p+q+1}T^*M\otimes \Vol^{k+l}),
$$
where $\Vol^k:=\begin{cases}(\Lambda^{n}TM)^{\otimes k}&\text{for
}k\geq 0\cr
(\Lambda^{n}T^*M)^{\otimes (-k)}&\text{for }k\geq
0.\end{cases}$

Observe that when we are not interested in rational
representations of the group of linear changes of coordinates but
allow ourselves to speak about infinitesimal transformations, we can
consider not only integer values of $k$ but any real or complex ones.

And the last (as we will see) invariant bilinear differential operator
$$
G: \Gamma(M, \Lambda^{p}T^*M\otimes \Vol^k)\times
\Gamma(M, \Lambda^{q}T^*M\otimes \Vol^l)\tto
\Gamma(M, \Lambda^{p+q-1}T^*M\otimes \Vol^{k+l}),
$$
was discovered in 1980, see \cite{G1}. The operator is a
generalization of two operators: the anti-symmetric Schouten
concomitant and its dual.

The same paper \cite{G1} contains the list of second and third order
differential operators. All of them are compositions of the exterior
differential $d$ and bilinear operators of orders $\leq 1$.

If $\dim M=1$, there exists one more new operator determined on the
weighted densities, not on the usual tensor fields, namely 
$$
T_{2}: f (dx)^{-2/3}, g (dx)^{-2/3}\mapsto
\left(2(f'''g-fg''')-3(f'g''-f''g')\right)(dx)^{5/3}.
$$
I discovered it in 1977, in my BS thesis. 

In 1979, Feigin and Fuchs
\cite{FF1} generalized it for the $m$-linear operators and in 1982
they classified all the multilinear anti-symmetric invariant
differential operators acting in the spaces of weighted densities on the
line \cite{FF2}.

The theorem on complete classification of differential operators
acting in the spaces of weighted densities on any manifolds is
announced in \cite{G1} and deposited to VINITI; here is a
slightly edited translation of the inaccessible deposition.

\ssec{2. Related results} Let an additional structure on
$M$ be fixed, e.g., a~volume, or a~symplectic structure, or a~contact structure.
One can consider differential operators invariant with respect to
transformations preserving the structure. We can, of course,
consider other types of structures, such as metrics or combinations of
several structures. But the method we use works well when the
Lie algebra of infinitesimal transformations is very asymmetric (has
more positive operators than negative ones) and close to simple.

For {\bf linear} (unary) operators, the complete classification was
obtained by Rudakov \cite{R1}, \cite{R2} for the general, volume preserving
and symplectic cases. I.~Kostrikin \cite{KoI} described the contact
case.

For {\bf bilinear} operators, the complete classification was obtained
for the general and volume preserving cases in \cite{G1} and for
symplectic case (partly) in \cite{G2}, \cite{G3}. I conjecture that
these partial results are final as far as indecomposable operators are
concerned, i.e., other, new, operators, if any, are compositions of the
ones already found. Observe that an explicit description of several
operators in symplectic case is to be given though their existence is proved \cite{G3}.

For the contact case, only small dimensions on supermanifolds
corresponding to some of the ``string theories'' are considered \cite{LKW}.

\ssec{3. Methods} 1) Reduction to canonical forms. For example, one can rectify any vector
field or a~volume form in a~vicinity of any non-singular point; in
other words, there are coordinates in which the components of these
tensor fields are constants. This means that the rational invariant
differential operators in the space $\Vect(M)$ or $\Vol(M)$ can only
be of order 0. In other words, they are algebraic, point-wise ones.

C.-L. Terng \cite{T} similarly proves that any rational invariant
differential operators in the spaces $C^\infty(M)$, or $\Omega^1(M)$,
or $\Omega^{n-1}(M)$ can be algebraically expressed via the exterior
derivative $d$.

Epstein \cite{E} similarly proved (making use of Cartan's results)
that --- and this is an important statement ---

{\sl any invariant differential operator on the quadratic forms can
be algebraically expressed via the curvature tensor and its covariant
derivatives}.

The tensor fields of more general form can not be reduced to a
canonical form by dimension considerations. Nevertheless, any tensor
field can be represented as a~sum of several tensor fields each of
which can be reduced to an {\it affine} form, i.e., to the form
in which the components of the field are vector-valued affine functions
$$
f(x)=a+\sum b_ix_{i}.
$$
It is precisely this fact that Palais \cite{P}, Leicher \cite{L} and
C.-L. Terng \cite{T} used to classify linear operators.

Kirillov in \cite{Ki1} uses another method. He considers any linear
invariant operator as a~morphism between two pairs of representations
of the Lie algebra of vector fields and its subalgebra $\fgl(n)$ of linear
vector fields. One further makes use of the sophisticated machinery
of representation theory, in particular, Laplace operators on finite
dimensional $\fgl(n)$-modules.

Here we come closer to the heart of the matter: the local problem
should be solved by local means and Lie algebras should replace
global discussions.

Rudakov \cite{R1} {\it started} with the infinitesimal problem. His
method applied to unary operators boils down to simple Linear Algebra only slightly seasoned with
some easy facts from representation theory and is applicable to
operators of any ``arity'', not only binary.

The method can be further applied to description of irreducible
representations of Lie algebras and superalgebras of vector fields. In
some cases the results directly follow from the description of
invariant linear differential operators and the Poincar\'e lemma or
its analogs (the general vector fields, see \cite{BL}). Kotchetkov
observed that sometimes (when no analog of Poincar\'e lemma holds) the situation is more subtle, see \cite{Ko}, \cite{Ko2}.

Bernstein showed \cite{BL} that local Rudakov's problem is equivalent to the
global one, initially formulated (somewhat vaguely) by Veblen and in
modern and lucid terms by Kirillov.

Let me describe Rudakov's method in more detail: it will be my main
tool in this paper.

\ssec{4. Rudakov's method for solution of Veblen's problem} Let $M$
be a~connected $n$-dimensional manifold over $\Ree$, and $\rho$ a
representation of $GL(n, \Ree)$ in a~finite-dimensional space $V$.
Denote by $T(\rho)$ or $T(V)$ the space of tensor fields of type
$\rho$ (or, which is the same, of type $V$), i.e., the collection of
the sections of the bundle over $M$ with fiber $V$ (over an open set
$U$). On $T(V)$, the group $\Diff M$ of diffeomorphisms of $M$ (the
local ones, which send $U$ into itself) naturally acts: let $J_A$ be
the Jacobi matrix of $A$ calculated in coordinates of points $m$ and
$A^{-1}(m)$, then set:
$$
A(t)(m)=\rho(J_A)(t(A^{-1}(m)))\quad \text{for } A\in\Diff M,\; m\in M,
t\in T(V). \eqno{(*)}
$$

Any operator $c: T(\rho_{1})\tto T(\rho_{2})$ is called {\it
invariant} if it commutes with the $\Diff M$-action.

It is instructive to compare Rudakov's and Kirillov's approaches to
Veblen's problem. First, they considered different categories, namely
Kirillov immediately confined himself to differential operators of
finite order and to tensor fields.

Rudakov allowed not only tensor fields but arbitrary jets and did not
bind the order of the (differential) operator. His result shows ({\it a
posteriori}) that

1) in spaces of jets higher than tensors (i.e., if the action depends
not only on first derivatives of the diffeomorphism, as in $(*)$, but
on higher derivatives) there are no invariant differential operators
(apart from scalar ones);

2) even if we consider arbitrary (irreducible) representations with
lowest weight vector, the restrictions on the weight that the
invariant operator requires for its existence imply that the
representation is finite-dimensional.

Observe that it is only due to the traditional reading of the term
``tensor field'' that we consider finite-dimensional representations.
It is more natural to consider, say, representations with vacuum vector
(the lowest for the tensor fields and the highest for the dual
spaces), though, strictly speaking, we have to consider indecomposable
representations in this infinite-dimensional setting.

Observe also that none of the researchers mention non-local invariant
operators: though we all know an example of such an operator --- the
integral --- it is unclear how to study them.

\ssec {5. The result} This paper contains

1) an enlarged reproduction of my Ph.D. thesis,
i.e., I give a~detailed proof of the classification of binary
differential operators listed in \cite{G1} and \cite{G2}.

2) The interpretation of some of the operators in terms of Lie
superalgebras seems to be new and might be of interest for theoretical
physicists.

Roughly speaking, the list of binary differential operators $D:
T(V_{1})\otimes T(V_{2})\to T(V_{3})$ invariant with respect to the
group of diffeomorphisms of $M$ runs as follows. Up to dualization and permutation of arguments,
the operators split into 9 types of order 1, four types of order 2 and
3 types of order 3. Operators of orders 2 and 3 are compositions of
1st order operators, except one indecomposable operator which only
exists for $n=1$. There are no operators of higher order.

Amazingly, almost all 1-st order operators determine a~Lie
superalgebra structure on their domain. Moreover, this Lie
superalgebra is almost simple: is a~central extension of a~simple one
or contains a~simple ideal of codimension 1.

3) In addition to the investigations from my thesis reproduced here,
I also considered the infinite-dimensional fibers. The result of
this consideration is discouraging: for 2-dimensional manifolds we do
not get ``really new'' operators (the operators we got earlier were
realized in functions polynomial fiber-wise; now we consider
nonpolynomial functions also but this is all); to consider manifolds of higher
dimensions seems to be a~wild problem.

\ssec {6. On open problems} A natural
generalization of the above Veblen--Rudakov's problem: consider
operators invariant with respect to other simple Lie algebras (or
superalgebras) of vector fields and consider operators of greater {\it
arity}: ternary, etc.

For the review of classification of unary operators (this task is
completely performed on manifolds and only
partly on supermanifolds), see \cite{L1}.

The case of binary operators is, so far, considered on symplectic
manifolds, see \cite{G2}, \cite{G3}, and on general and certain
contact supermanifolds \cite{LKW}. Both results are partial.

The exceptional bilinear operator of order 3 was generalized in
\cite{FF1}, \cite{FF2}, where $m$-ary anti-symmetric operators on the line
are classified.

The generalization of the problem in all the directions mentioned is
desirable, but we give the priority to the operators invariant with
respect to the Lie algebra that preserves the contact form on
manifolds and various structures on the supercircle as having
more immediate applications.

\ssec{A mysterious operator} Since in the absence of even coordinates,
all operators in superspaces of tensors or jets (with finite
dimensional fiber) are {\it differential} ones (even the integral),
Leites hoped that having worked out this finite-dimensional model one
could, by analogy, find new non-local invariant operators for manifolds
as well. The arity of such operators is, clearly, $>1$. So far, no
such operator is explicitly written except an example of a~{\it symbol of such a
binary} operator acting in the spaces of certain tensors on the line;
see Kirillov's review \cite{Ki3}.

Recently, this operator was demystified, see \cite{IoMa}. Later on, Bouarroudj and Leites classified bilinear differential operators on 1-dimensional supervariety over algebraically closed fields of characteristic $p>0$ and found several generalizations of what they called, after Feigin and Fuchs, the ``Grozman operator", see \cite{BoLe}.

\section{The list of operators}

Let $\Omega^i=T(\Lambda ^i (\id))$ be the space of differential
$i$-forms, $\Omega^{\bcdot}=\mathop{\oplus}\Omega^i$. Recall that
$n= \dim M$; let $\Vol=\Omega^n$ and let $\Vol^\lambda$ for
$\lambda\in \Cee$ be the space of $\lambda$-densities. This is a
rank 1 module over functions $\cF=\Omega^0$ with generator
$\vvol_{x}^\lambda$. Observe that the action of $\Diff M$ is not
defined on $\Vol^\lambda$ unless $\lambda$ is integer, but the Lie
algebra $\fvect(n)$ naturally acts on $\Vol^\lambda$ for any
$\lambda$ by the formula
$$
L_{D}(\vvol_{x}^\lambda)=\lambda\Div (D)\vvol_{x}^\lambda\text{ for
any }D\in\fvect(M).
$$
We will consider this later wider problem: classification of
$\fvect(M)$-invariant differential operators.

\ssec{1.0. Zero order operators} Obviously any zero order
differential operator
$$
Z: T(V_{1})\otimes T(V_{2})\tto T(W)
$$
is just a~scalar one and is the uniquely defined extension of a~morphism of $\fg_{0}$-modules in $\Hom(V_{1}\otimes V_{2}, W)$.

\ssec{1.1. First order operators}
$$
P_{1}: \Omega^r\otimes T(\rho_2) \tto T(\rho_3)\quad (w,t) \mapsto Z(dw,
t);\eqno{(P_{1})}
$$
where $Z$ is the zero-th order operator --- extension of the
projection $\rho_1\otimes \rho_2\tto \rho_3$ onto the irreducible
component; the operator $P_{1}^{*2}$ is of the same form, whereas
$P_{1}^{*1}$ is of the form
$$

\begin{array}{l}
P_{1}^{*1}: T(\rho_1)\otimes T(\rho_2) \tto \Omega^r \qquad
P_{1}^{*1}(t_1,t_2) \mapsto d(Z(t_1,t_2)),\\
 \text{where } Z:
T(\rho_1)\otimes T(\rho_2) \tto \Omega^{r-1}\text{ for } r>0.
\end{array}
$$
The Lie derivative:
$$
P_2: \Vect \otimes T(\rho) \tto T(\rho); \eqno{(P_{2})}
$$
the operator $P_{2}^{*2}$ is also the Lie derivative, whereas
$P_{2}^{*1}$ is Schouten's ``lagrangian concomitant''.

Schouten's ``symmetric concomitant'' or the Poisson bracket:
$$
P_3 = PB: T(S^p(\id^*)), T(S^q(\id^*)) \tto
T(S^{p + q -1}(\id^*)); \eqno{(P_{3})}
$$
for $p=1$ the operator reduces to the Lie derivative, for $p=1$ to
$P_1$; the duals of $P_3$ are also of the same (up to the permutation of arguments $T$) form:
$$
P_3^{*2} : T(S^p(\id^*)), T(S^q(\id))\otimes\Vol^*\tto
T(S^{p + q -1}(\id))\otimes\Vol^*.
$$
The {\it Nijenhuis bracket}. This bracket is a~linear combination of
operators $P_1$, $P_1^*1$, their composition with the permutation of arguments (a.k.a. twist)
operator 
\[
T: T(V)\otimes T(W) \tto T(W)\otimes T(V), 
\]
and an
``irreducible'' operator sometimes denoted in what follows by $N$
$$
P_4: \Omega ^p\otimes_c\Vect, \Omega ^q\otimes_c\Vect \tto
\Omega ^{p + q}\otimes_c\Vect \eqno{(P_{4})}
$$
defined as follows
$$
\begin{array}{l}
 P_4(\omega_1\otimes D_1, \omega _2\otimes D_2)=(\omega_1\wedge \omega
 _2)\otimes [D_1, D_2]+\\
\left (\omega_1\wedge L_{D_{1}}(\omega _2)+(-1)^{p(\omega_1)}d\omega_1\wedge
\iota_{D_{1}}(\omega _2)\right )\otimes D_2+\\
\left (-L_{D_{2}}(\omega_1)\wedge \omega
_2+(-1)^{p(\omega_2)}\iota_{D_{2}}(\omega _1)d\wedge \omega_2
\right )\otimes D_1.
\end{array}
$$
The invariance of the Nijenhuis bracket is a~corollary of the
following observation. It is evident that for a~fixed
$\omega_1\otimes D_1\in\Omega^k\otimes_c\Vect$, the operator
$D':\Omega^k\otimes_c\Vect \times \Omega^{\bcdot}\tto \Omega^\bcdot$
given by the formula (here $\overline B(a, b):= B(b, a)$ is the twisted
operator, $\iota_D$ is the inner derivation along $D$)
$$
\begin{array}{l}
D'(\omega_1\otimes D_1, \omega _2)=D^{*1}(\omega_1\otimes D_1, \omega _2)+
\overline D(\omega_1\otimes D_1, \omega _2)=\\
d(\omega_1\wedge L_{D_{1}}(\omega _2))+(-1)^{p(\omega_1)}\omega_1\wedge
\iota_{D_{1}}(d\omega _2)=\\
d\omega_1\wedge \iota_{D_{1}}\omega
_2+(-1)^{p(\omega_1)}\omega _1\wedge L_{D_{1}}(\omega _2).
\end{array}
$$
is a~superdifferentiation of the supercommutative superalgebra $\Omega^{\bcdot}$.

Observe that the $\fgl(n)$-module $\Lambda^k(\id)\otimes \id^*$ is
reducible:
$$
\Lambda^k(\id)\otimes \id^*=R(1, 0, \dots, 0, \underbrace{-1, \dots,
-1}_{k})\oplus R(0, \dots, 0, \underbrace{-1, \dots, -1}_{k-1}).
$$
Therefore, the operator $N$ splits into the direct sum of several
operators. One of them, that we did not consider before, will be
denoted $N$ or $P_4$, namely the projection onto the first component:
$$
\begin{array}{l}
 N: T(R(1, 0, \dots, 0, \underbrace{-1, \dots,
-1}_{k}))\times T(R(1, 0, \dots, 0, \underbrace{-1, \dots,
-1}_{l}))\tto \\
T(R(1, 0, \dots, 0, \underbrace{-1, \dots,
-1}_{k+l})).
\end{array}
$$

There is also a~dual operator:
$$
\begin{array}{l}
 N^{*2}: T(R(1, 0, \dots, 0, \underbrace{-1, \dots,
-1}_{k}))\times T(R(0, \dots, 0, \underbrace{-1, \dots,
-1}_{l}), -2)\tto \\
T(R(0, \dots, 0, \underbrace{-1, \dots,
-1}_{k+l}), -2).\end{array}
$$

The following operator is just a~composition of the exterior
derivative and a~zero order operator:
$$
P_5: \Omega ^p, \Omega ^q \tto \Omega ^{p + q + 1}; \quad
\omega_1, \omega_2\mapsto (-1)^{p(\omega_1)}a(d\omega_1\omega_2) +
b(\omega_1d\omega_2), \text{ where } a, b\in\Cee. \eqno{(P_5)}
$$

Let $|\mu|^2 + |\nu|^2 \neq 0$. Define
$$
P_6^\Omega : \Omega ^p_{\mu}, \Omega ^q_{\nu}\tto \Omega
^{p+q+1}_{\mu+\nu}
$$
by setting
$$
\omega_1\vvol^\mu, \omega_2\vvol^\nu\mapsto \left (\nu
(-1)^{p(\omega_1)}d\omega_1\omega_2 -
\mu\omega_1d\omega_2\right)\vvol^{\mu+\nu}. \eqno{(P_6)}
$$
Denote the Schouten bracket:
$$
P_7 : L^p, L^q \tto L^{p + q - 1}.\eqno{(P_7)}
$$

Define a~generalization $P_8: L^p_{\mu}, L^q_{\nu} \tto
L^{p+q-1}_{\mu+\nu}$ of the Schouten bracket (on manifolds, for $p
+ q \leq n$; on supermanifolds of dimension $n|1$, for $p,
q\in\Cee$) by the formula
$$

\begin{array}{l}
 X\vvol^\mu, Y\vvol^\nu \mapsto \left ((\nu-1)(\mu + \nu -1) \Div X \cdot Y +
(-1)^{p(X)}(\mu-1)(\mu + \nu -1)X \Div Y -\right .\\
\left . (\mu- 1)(\nu -1)
\Div(XY)\right )\vvol^{\mu+ \nu},\end{array}\eqno{(P_8)}
$$
where the {\it divergence of a~polyvector field} is best described in
local coordinates $(x, \check x)$ on the supermanifold $\check M$
associated to any manifold $M$, cf. \cite{BL}.

The operators dual to $P_6$, $P_7$, $P_8$ are, as is not difficult to
see, of the same form, respectively.

\ssec{1.2. Operators of order $>1$} All of them are reduced to
compositions of operators of orders $\leq 1$:
$$
S_{1}: \Omega ^p\times\Omega ^q\tto T(R(0, \dots, 0, \underbrace{-1, \dots,
-1}_{k}, \underbrace{-2, \dots , -2}_{l}), \text{ where }
k+2l=p+q+2, \eqno{(S_{1})}
$$
defined to be
$$
S_{1}(\omega_1, \omega_2)=Z(d\omega_1, d\omega_2);
$$
$$
S_{1}^{*}: \Omega ^p\times T(R(1, \dots, 1, 0, \dots, 0, -1, \dots,
-1)\tto \Omega ^q, \eqno{(S_{1}^{*})}
$$
$$
S_{1}^{*}(\omega, t)=dZ(\omega, t);
$$
$$
S_{2}: \Omega ^{n-1}\times\Omega ^p\otimes \Vol^k\tto \Omega ^{p+1}\otimes
\Vol^{k+1}, \eqno{(S_{2})}
$$
defined to be
$$
S_{2}(\omega, t)=F(d\omega, t);
$$
$$
S_{2}^{*}: \Omega ^{p}\otimes \Vol^k\times\Omega ^{n-1-p}\otimes
\Vol^{-k-1}\tto \Omega ^{1}, \eqno{(S_{2}^{*})}
$$
defined to be
$$
S_{2}^{*}(a, b)=dF(a, b);
$$
$$
T_{1} \Omega ^{n-1}\times\Omega ^{n-1}\tto \Omega ^{1}\otimes
\Vol^{2}, \eqno{(T_{1})}
$$
defined to be
$$
T_{1}(\omega_1, \omega_2)=F(d\omega_1, d\omega_2);
$$
$$
T_{1}^{*}: \Omega ^{n-1}\times\Omega ^{n-1}\otimes
\Vol^{-2}\tto \Omega ^{1}, \eqno{(T_{1}^{*})}
$$
defined to be
$$
T_{1}^{*}(\omega, t)=dF(\omega, t);
$$

If we abandon requirement of rationality of densities in the
definition of operators $F$, $G$ and~ $S_{2}$, then for $n=1$ we
obtain one more (irreducible, i.e., not factorizable in a~composition)
invariant operator
$$
T_{2}: \Vol^{-2/3}\times\Vol^{-2/3}\tto
\Vol^{5/3}, \eqno{(T_{1})}
$$
defined to be
$$
T_{2}:(f\vvol^{-2/3}, g\vvol^{-2/3})\longmapsto (2f'''g-2fg'''+3f''g'-3f'g'')\vvol^{5/3}.
$$

\SSSbegin{1.2.1}{Theorem}\label{Th} Every bilinear invariant differential operator acting
in tensor fields on a~connected smooth manifold is a~linear
combination of the above operators and the ones obtains from them
by a~transposition of the arguments. \end{Theorem}

\section{The beginning of the proof}

Consider the group $G=\Diff M$ of local diffeomorphisms of $M$.
In some sense, its Lie algebra is $\fg=\fvect (M)$, the Lie
algebra of vector fields on $M$. Let $G(x_0)$ be the stabilizer
of $x_0\in M$; its Lie algebra is $\fg(x_0)=\{\xi\in
g\mid\xi(x_0)=0\}$. There exists a~neighborhood $U$ of $x_0$ over
each point of which the fibers can be identified with a
``standard" fiber $V$, then $\fg$ acts on the tensor fields, the
elements from $T(V)$, via the formula
$$
L_{\xi}(\varphi\otimes v)=\sum_{i=1}^{n}\xi_{i}\frac{\partial\varphi}
{\partial x_i}\otimes v+\sum^n_{i, j=1}\frac{\partial\xi_i}{\partial
x_j}\otimes \rho(E^j_i)v, \eqno {(1)}
$$
where $\varphi\in C^{\infty}(U)$, $v\in V$,
$\xi=\sum^n_{i=1}\xi_i\frac{\partial}{\partial x_i}\in \fg$ and
$\{E^j_i\}^n_{i, j=1}$ is the standard basis of $\fgl(n)$
consisting of matrix units. The space
$I(V^*)=\Kee(\partial_1, \dots, \partial_m)\otimes V^*$ of differential
operators whose coefficients are linear functionals on $V$ is a~$\fg$-invariant
subspace of $(T(V))^*$. The pairing of $I(V^*)$ with $T(V)$ is determined by the
formula
$$
\langle\partial\otimes v', f\otimes v\rangle=\partial f\mid_{x=0}\langle v', v
\rangle,
$$
where $\partial\in\Kee[\partial_1, \dots, \partial_m]$, $f\in
C^{\infty}(M)$, $v\in V$, $v'\in V^*$ and $f\otimes v$ is a
representation of the section $s\in T(V)$ in coordinates. The action
of $\fg$ is found from the formula
$$

\begin{array}{l}
 \langle I_{\xi}(\partial\otimes v'), f\otimes v\rangle=
-\langle\partial\otimes v', L_{\xi}(f\otimes v)\rangle=\\
-(\partial\circ\xi f)\mid_{x=0}\langle v', v\rangle
-\sum_{i, j}\partial(\frac{\partial\xi_i}{\partial x_j} f)\mid_{x=0}
\langle v', \rho(\varepsilon ^j_i)v\rangle=\\
\langle-\partial\circ\xi\mid_{x=0}\otimes v'+\sum_{i, j}
(\partial\circ\frac{\partial\xi_i}{\partial x_j})\mid_{x=0}\otimes\rho^*
(\varepsilon ^j_i)v', f\otimes v\rangle.
\end{array}
$$
Thus,
$$
I_{\xi}(\partial\otimes v')=(\partial\circ\xi)\mid_{x=0}\otimes
\sum_{i, j}(\partial\circ\frac{\partial\xi_i}{\partial x_j})\mid_{x=0}
\otimes\rho^*(E^j_i)v'.\eqno{(2)}
$$
The space $I(V^*)$ is graded $I(V^*)=\otimes_{m\geq 0}I^m(V^*)$, where
$I^m(V^*)$ consists of homogeneous polynomials
of degree $m$ in $\partial_1, \dots, \partial_n$. Let
$$
\cI_m(V^*)=\oplus^{m}_{k=0}I^k(V^*)
$$
denote the space of polynomials of degree $\leq m$. Observe, that each
$\cI_m(V^*)$ is $\fg(x_0)$-invariant, in particular, $\cI_0(V^*)=V^*$.

Together with $\fg(x_0)$, consider the Lie algebra $\cL_0=\fvect(n)$ of polynomial
vector fields on $\Kee^n$
that vanish at the origin, $x_0$. Determine the $\cL_0$-action on $I(V^*)$ by the
same formula (2). Clearly, with $\cL_0$ a~grading is associated
$$
\cL_0=L_0\oplus L_1\oplus\dots,
$$
where $L_m$ consists of vector fields whose coefficients (i.e., the
coefficients of the $\partial_i$) are homogeneous polynomials of
degree $m+1$. The Lie subalgebra $L_0$ is isomorphic to $\fgl(n)$
under the correspondence
$$
x_i\partial_j\longleftrightarrow E^i_j.
$$
Observe that if $\xi\in L_0$, then $I_\xi(1\otimes v)=1\otimes\rho^*(\xi)v$,
whereas $\xi\in\cL_1=L_1\oplus L_2\oplus\dots$ annihilates $\cI_0(V^*)$.

\begin{Lemma}\label{L2.1} $\cI_m(V^*)$ is the annihilator of
$Z_{x_0}^{m+1}(V)=I^{m+1}_{x_0}\cdot T(V)$.
\end{Lemma}

\begin{proof} If $\deg\partial\leq m$ and $\varphi\in I^{m+1}_{x_0}$, then
$\partial(\varphi s)(0)=0$. Hence, $\langle\cI_m(V^*),
Z_{x_0}^{m+1}(V)\rangle=0$,
but $\dim\cI_m(V^*)=\codim Z^{m+1}_{x_0}(V)$; hence, $\cI_m(V^*)$ is the
annihilator of $Z^{m+1}_{x_0}(V)$.

Let $D: T(V_1)\to T(V_2)$ be a~$G$-invariant differential operator of order
$m$. Then, 
$$
D(Z^{m+1}_{x_0}(V_1))\subset Z^1_{x_0}(V_2),
$$
and therefore 
$D^*(\cI_0(V^*_2))\subset\cI_m(V^*_1)$.
\end{proof}

\begin{Remark} By means of the standard theorems of Linear Algebra one
can prove that for any $\fg(x_0)$-invariant operator $D^0:
V^*_2\to\cI_m(V^*_1)$, there exists a~unique $\fg$-invariant operator
${D: T(V_1)\to T(V_2)}$ such that $D^*\mid_{\cI_0(V^*_2)}=D^0$.
The $\fg(x_0)$-action on $\cI_m$ depends only on the first $m+1$
derivatives (from the 0-th to $(m+1)$-st inclusively) at the origin
$x_0$ of
the vector fields from $\fg(x_0)$. Hence, the sets of operators
$$
\{I_{\xi}:\xi\in g(x_0\} \; \text{ and }\; \{I_{\xi}:\xi\in\cL\}
$$
coincide, and therefore the $\fg(x_0)$-invariance of the operator
$D^*: V_2\to\cI_m(V^*_1)$ is
equivalent to its $\cL$-invariance.
\end{Remark}

The fact that $D^*(V_2^*)=W\subset \cI_m(V_1^*)$ is an
$\cL$-submodule isomorphic to $V^*_2$ (which is clear, since $V^*_2$
is irreducible) implies that

a) $W$ is an $L_0$-submodule isomorphic to $V^*_2$;

b) $\cL_1$ annihilates $W$.

The problem of description of $\fg$-invariant differential operators is, therefore,
equivalent to the following problem:

{\sl in $I(V^*_1)$, find all $L_0$-submodules
isomorphic to $V^*_2$ and such that $\cL_1$ annihilates them.}

The vectors that $\cL_1$ annihilates will be called {\it singular}
ones. Thus, our problem is to describe highest weight singular vectors.

In case of the bilinear operators $B: T(V_1)\otimes T(V_2)\tto T(V_3)$
the above procedure should be modified as follows. Observe that
$T(V_1)\otimes T(V_2)$ is a~$C^{\infty}(M)\otimes
C^{\infty}(M)$-module, whereas $\cI_{\rho_0}=I_{\rho_0}\otimes 1+\otimes
I_{\rho_0}$ is a~maximal ideal of $C^{\infty}(M)\otimes
C^{\infty}(M)$. Let
$$
Z^m_{x_0}(V_1, V_2)=\cI^m_{x_0}(T(V_1)\otimes T(V_2)).
$$
The space
$$

\begin{array}{l}
I(V^*_1, V^*_2)=I(V^*_1)\otimes I(V^*_2)=\Kee[\partial'_1, \dots, \partial'_n]\otimes
V^*_1\otimes\Kee[\partial''_1, \dots, \partial''_n]\otimes V^*_2=\\
\Kee[\partial'_1, \dots, \partial''_n]\otimes(V^*_1\otimes V^*_2)
\end{array}
$$
is graded by the total degree of the polynomials in $\partial'$ and
$\partial''$. Clearly,
$$
I(V^*_1, V^*_2)=\oplus^{\infty}_{m=0}I^m(V^*_1, V^*_2);
$$
and
$$
\cI_m(V^*_1, V^*_2)=\oplus^{m}_{i=0}I^i(V^*_1, V^*_2)
$$
is the annihilator of $Z^{m+1}_{x_0}(V_1, V_2)$.
If $B$ is a~$\fg$-invariant bilinear operator of order $m$, then
$$
B(Z^{m+1}_{x_0}(V_1, V_2))\subset Z^1_{x_0}(V_3),
$$
hence,
$$
B^*(\cI_0(V^*_3))\subset\cI_m(V^*_1, V^*_2).
$$
Therefore, to find all such $B$'s it remains to find in $I(V^*_1,
V^*_2)$ all $L_0$-submodules annihilated by~ $\cL_1$.

\section{Solution $(n=1)$}

All irreducible finite-dimensional and diagonalizable modules of
$\fgl(1)$ are 1-dimensional; let $V^*_1$ and $V^*_2$ be such
modules. Let $v\in V^*_1$ and $w\in V^*_2$ be nonzero vectors of
weight $l$ and $m$, respectively, i.e.,
$$
(x\partial)v=lv, \quad (x\partial)w=mw.
$$
Since for $n=1$, there is no notion of highest weight $\fgl(1)$-vector,
it suffices to describe the singular vectors in $I(V^*_1, V^*_2)$.
Since $\cL_1$ is generated by $\varepsilon _1=x^2\partial$ and
$\varepsilon _2=x^3\partial$, it suffices to find all weight solutions
of the system
$$
\varepsilon _1f=0, \quad \varepsilon _2f=0
$$
for a~homogeneous vector $f\in I(V^*_1, V^*_2)$ of
degree $d$:
$$
f=\sum_{i+j=d}\frac{1}{i!j!}c_i\partial^iv\otimes\partial^jw,
$$
where the factor $\frac{1}{i!j!}$ is inserted for further convenience.

Observe that
$$
(x^2\partial)(\partial^i\otimes v)=2i\partial^{i-1}\otimes(x\partial)v-i(i-1)
\partial^{i-1}\otimes v=i(2l-i+1)\partial^{i-1}v;
$$
$$

\begin{array}{l}
 (x^3\partial)(\partial^iv)=3i(i-1)\partial^{i-2}\otimes(x\partial)v-i(i-1)
(i-2)\partial^{i-2}\otimes v=\\
i(i-1)(3l-i+2)\partial^{i-2}v.
\end{array}
$$
Hence,
$$

\begin{array}{l}
0=(x^2\partial)f=\sum_{i+j=d}\frac{c_i}{i!j!}(i(2l-i+1)\partial^{i-1}v\otimes
\partial^jw+j(2m-j+1)\partial^iv\otimes\partial^{j-1}w)=\\
\sum_{i+j=d}\left(\frac{c_{i+1}}{i!j!}(2l-i)\partial^iv\otimes\partial^{j-1}w+\frac{c_i}{i!j!}
(2m-j+1)\partial^iv\otimes\partial^{i-1}w\right )
\end{array}
$$
implying
$$
(2l-i)c_{i+1}+(2m-j+1)c_i=0
$$
and
$$

\begin{array}{l}
 0=(x^3\partial)f=\\
 \sum_{i+j=d}\frac{c_i}{i!j!}(i(i-1)(3l-i+2)\partial^{i-2}v
\otimes\partial^jw+j(j-1)(3m-j+2)\partial^iv\otimes\partial^{j-2}w=\\
\sum_{i+j=d}\frac{1}{i!j!}\left (c_{i+2}(3l-i)\partial^iv\otimes\partial^{j-2}w+
c_i(3m-j+2)\partial^iv\otimes\partial^{j-2}w\right )
\end{array}
$$
implying
$$
(3l-i)c_{i+2}+(2m-j+2)c_i=0.
$$
Let $i+j=d$, $1\leq i$, $j\leq d-1$. The above formulas impose the
constraints on $c_{i-1}$, $c_i$
and $c_{i+1}$:
$$
\left\{\begin{matrix}(2l-i)c_{i+1}+(2m-j+1)c_i=0\cr
 (2l-i+1)c_i+(2m-j)c_{i-1}=0\cr
(3l-i+1)c_{i+1}+(3m-j+1)c_{i-1}=0.\cr
\end{matrix}\right.
$$
The determinant of this system is
$$

\begin{array}{l}
 \triangle_i=\left|\,
\begin{matrix}2l-i&2m-j+1&0\cr
0& 2l-i+1& 2m-j\cr
3l-i+1& 0& 3m-j+1\cr
\end{matrix}\right |=\\
(2l-i)(2l-i+1)(3m-j+1)+(2m-j)(2m-j+1)
(3l-i+1).
\end{array}
$$

Observe that if $c_i=0$ and $c_{i+1}=0$, then
\[
(2m-j)c_{i-1}=(3m-j+1)c_{i-1}=0
\]
 but since $2m-j$ and $3m-j+1$ cannot
vanish simultaneously for $j\geq 0$, it follows that $c_{i-1}=0$.
Hence, if two neighboring coefficients in a~row vanish, then the
coefficients neighboring them also vanish.

So, for a~nonzero solution to exist, it is necessary (but not
sufficient) that $\triangle_i=0$ for any $i$ such that $1\leq i\leq d-1$.

Set $x=2l-i+1$ and $y=2m-j+1$. Then, in terms of $\cD:=x+y=2l+2m-d+2$
we have
$$

\begin{array}{l}
 \triangle_i=\left|\,
 \begin{matrix}
 x-1& y& 0\cr
 0& x& y-1\cr
 x+l& 0& y+m\cr
\end{matrix}\right |=
(x+l)y(y-1)+(x+1)x(y+m)=\\
x^2y+xy^2-2xy+ly(y-1)+mx(x-1)=\\
x^2(\cD-x)+x(\cD-x)^2-2x(\cD-x)+l(\cD-x)(\cD-x-1)+mx(x-1)=\\
x^2(l+m+2-\cD)+x(\cD-2\lambda\cD+l-m)+l\cD^2-l\cD=0. 
\end{array}
$$

If $d\geq 4$, then the quadratic equation has $\geq 3$ solutions which
is only possible if all the coefficients vanish:
$$
\text{$\lambda +\mu+2-\cD=0$ or $l+m+2=2l+2m-d+2$}
$$
implying $d=l+m$ and
$$
\cD^2+2\lambda\cD-2\cD+l-m=-l^2+m^2-l+m=(m-l)(l+m+1)=0.
$$

But $m+l=d-1\neq 0$ yields $l=m$; so
\[
l\cD^2-l\cD=l\cD(\cD-1),
\]
but $\cD=l+m+2=d+2\geq 6$; hence, $l=0$, but then
$m=0$ as well, implying $d=l+m=0$. This is a~contradiction.

There are no nonzero solutions, hence, on the 1-dimensional manifold,
there are no bilinear operators of order $>3$.

For $d=3$, the equation $\triangle_i=0$ has two roots, 1 and 2:
$$
\left\{\begin{matrix}
(2l-1)\cdot 2l\cdot (3m-1)+(2m-2)(2m-1)\cdot 3l=0\cr
(2l-2)(2l-1)\cdot 3m+(2m-1)\cdot 2m(3l-1)=0.\cr\end{matrix}\right .
$$

a) $l=0$. Then, $6m=4m^2+2m$ implying either $m=0$ or $m=2$.
Thus, there are solutions $(0, 0)$, $(0, 2)$ and a~symmetric solution $(2, 0)$.

b) $l\neq0$, $m\neq0$. Then, 
$$
\begin{matrix}(2l-1)(3m-1)+3(m-1)(2m-1)=0\Longrightarrow \\
3(2l-1)(3m-1)(l-1)+9(m-1) (2m-1)(l-1)=0;\\
3(l-1)(2l-1)+(2m-1)(3l-1)=0\Longrightarrow\\
3(l-1)(2l-1)(3m-1)+(2m-1)(3l-1) (3m-1)=0;
\end{matrix}
$$
implying
$$
9(m-1)(l-1)(2m-1)=(3m-1)(2m-1)(3l-1).
$$
Hence, either $m=\frac{1}{2}$ and then $(2l-1)(3m-1)=0\Longrightarrow
l=\frac{1}{2}$ or

$9(m-1)(l-1)=(3m-1)(3l-1)$, i.e., $9lm-9l-9m+9=9lm-3l-3m+1\Longrightarrow l+m=
\frac{4}{3}$;

$3(l-1)(2l-1)=(\frac{8}{3}-2l-1)(3l-1)$,
$2l=\frac{4}{3}$ so, finally, $l=m=\frac{2}{3}$.

But the condition was not a~sufficient one; the complete condition is
$$
\left\{\begin{matrix}
(2m-2)c_0+2lc_1=0, \cr
(2m-1)c_1+(2l-1)c_2=0, \cr
2mc_2+(2l-2)c_3=0, \cr
(3m-1)c_0+3lc_2=0, \cr
3mc_1+(3l-1)c_3=0.\cr\end{matrix}\right .
$$
It is routine to verify that in all the cases except for $l=m=\frac{1}{2}$
there is a~solution and this solution is unique up to multiplication
by a~constant; whereas for $l=m=\frac{1}{2}$ there are no
solutions. (This is in agreement with our list of operators.)

Let $d=2$, then
the equation $\triangle_1=0$ is equivalent to
$$
(2l-1)\cdot 2l\cdot 2m+(2m-1)\cdot 2m\cdot 3l=0.
$$

a) $l=0$, $m$ is arbitrary. The matrix $\left(\begin{matrix}
-1&2m&0\cr
0&0&2m-1\cr
0&0&3m\cr
\end{matrix}\right)$ is of rank 2, hence, there
is one solution in this case (namely, the operator $S_2$, which
turns into $S_1$ for $m=0$).

The case $m=0$ is similar.

b) $l\neq 0$, $m\neq 0$. Then, 
$(2l-1)+(2m-1)=0$, $l+m=1$, $m=1-l$ and the rank of
$\left(\begin{matrix}2l-1&2-2l&0\cr
0&2l&1-2l\cr
3l&0&3-3l\cr
\end{matrix}\right)$ is always
equal to 2; hence, we have only one
operator, $S^{*1}_2$.

In case $d=1$, there remains one condition $2l\cdot c_0+2m\cdot c_1=0$.
If $l$ and $m$ do not vanish simultaneously, we have just one operator,
$P_4$, and if $l=m=0$, then we have two operators (both of type $P_1$):
$$
B(\varphi, \psi)=a\varphi d\psi+bd\varphi\cdot \psi.
$$

\section{Solution $(n=2)$}
We denote the operator $\xi\in\cL$ acting on
$I(V_1^*)\otimes I(V^*_2)$ by the same symbol $\xi$ as the element
itself, but $\xi'$ indicates that $\xi$ acts on the first factor of
$I(V_1^*)\otimes I(V_2^*)$ whereas $\xi''$ acts
only on the second factor.

We identify the linear vector fields with $2\times 2$
matrices and use the following shorthand notations:
$$
X_+=x_1\partial_2, \quad X_-=x_2\partial_1, \quad h_1=x_1\partial_1, \quad
h_2=x_2\partial_2.
$$
The weights of representations $\rho_1^*$, $\rho_2^*$ and $\rho_3^*$ with
respect to $h_1$ and $h_2$ will be denoted by
$$
\bar\lambda=(l_1, l_2), \quad \bar\mu=(m_1, m_2), \quad \bar\nu=(n_1, n_2);
$$
the weights with respect
to $h_1-h_2$ (in other words, with respect to $\fsl(2)$) are:
$$
\lambda=l_1-l_2, \quad \mu=m_1-m_2, \quad \nu=n_1-n_2.
$$
Sometimes the subscript 2 will be omitted.

In $V^*_1$ and $V^*_2$, fix weight bases $v_0, v_1, \dots, v_{\lambda}, \dots$ and
$w_0, w_1, \dots, w_\mu, \dots$ such that
$$
X_+v_i=(\lambda-i+1)v_{i-1}\text{ for $i<0$}, \quad X_+v_0=0
$$
$$
X_-v_i=(i+1)v_{i+1}\begin{cases}\text{for all }i&
\text{if }\lambda\not\in\Zee_+,\cr
\text{for }i<\lambda&\text{if }\lambda\in\Zee_+.\cr
\end{cases}
$$
If $\lambda\in\Zee_+$, we set $X_-v_\lambda=0$; moreover, we only
consider $v_0, \dots, v_\lambda$. Similar formulas apply to the
$w_0, \ldots\in V_2^*$.

In $I(V_1^*, V_2^*)$, the weight basis consists of
vectors
$$
(\partial'_1)^\alpha(\partial'_2)^\beta(\partial''_1)^\gamma
(\partial''_2)^\delta v_i\otimes w
$$ for the above $v_i, w_j$. Observe that the
weights of $\partial'_1$ and $\partial''_1$ are $(-1, 0)$, the
weights of $\partial'_2$ and $\partial''_2$ are $(0, -1)$.
The weight vector $f$ is a~highest one if $X_+f=0$.

The irreducible finite-dimensional $\fgl(2)$-module is determined by
its highest weight vector up to an isomorphism. In particular, in order
to find all finite-dimensional irreducible $L_0$ submodules of
$I(V_1^*, V_2^*)$ it suffices to find all the highest weight vectors
in $I(V_1^*, V_2^*)$.

Further, among these vectors, we have to find the singular vectors,
the ones that $\cL_1$
annihilates. To this end, it suffices to solve the system of equations
$(x_2^2\partial_1)f=0$ and $(x_2^2\partial_2)f=0$ because every
element from $\cL_1$ can be expressed in terms of $x_2^2\partial_1$,
$x_2^2\partial_2$ and $X_+$.

Thus, our problem is to find the homogeneous (with respect to weight) solutions of
the system
$$
X_+f=(x_2^2\partial_1)f=(x_2^2\partial_2)f=0.\eqno {(3)}
$$
The action of $\cL$ on $I(V_1^*, V_2^*)$ is compatible with the
grading, and therefore any solution of system is the sum of
homogeneous solutions. We will seek only homogeneous solutions.
Observe that the homogeneity degree of the singular vector
coincides with the order of the corresponding differential
operator.

We will look for solutions in the form
$$
f=\sum P_i(\partial'_1, \partial'_2, \partial''_1, \partial''_2)\otimes u_i,
$$
where
$P_i$ are monomials of degree $d$ and $u_i\in V_1^*\otimes V_2^*$
(in what follows we will often omit the sign of the tensor product).
                   
\begin{Lemma}\label{L4.1} Weight solutions of the equation $(X_+)^{d+1}u=0$, where
$u\in V_1^*\otimes V_2^*$, are of the form
$$
u=\sum^{s}_{i=0}(-1)^i(\lambda-i)!(\mu-s+i)!P(i)v_i\otimes w_{s-i},
\eqno {(4)}
$$
where $P(i)$ is a~polynomial of degree not greater than $d$.
\end{Lemma}

I will denote the elements $u$ of the form (4) by $(s; P(i))$. The
weight of such an element is equal to $(l_1+m_1-s, l_2+m_2+s)$.

Proof of the lemma follows from the formula
$$
X_+(s, P(i))=(s-1, P(i)-P(i+1))
$$
(that is the $\deg (P(i)-P(i+1))-\deg P(i) =1$).
Observe also that the action of $X_-$ is as follows:
$$
x'_-(s, P(i))=(s+1, i(i-\lambda-1)P(i-1))\qquad (i=0, 1, \dots,
s+1).\eqno{(5)}
$$
$$
x''_-(s, P(i))=(s+1, -(i+\mu-s)(i-s-1)P(i))\qquad (i=0, 1, \dots, s+1).\eqno{(6)}
$$
In what follows we will see
that the elements $u_j$ which correspond to monomials of the form
$P(\partial'_2, \partial''_2)$ of
degree $d$ in
the decomposition (4) of $f$ satisfy
$$
(X_+)^{d+1}u=0.
$$
The weight of
such $u$ is equal to $(l_1+m_1-s, l_2+m_2+s)$, and
therefore the weight of $f$ is equal to $(l_1+m_1-s, l_2+m_2+s-d)$.

If $B: (T(V_1), T(V_2))\to T(V_3)$ is a~homogeneous differential
operator of order $d$, then its first and second duals are of
the same order. This can be deduced by integrating by parts
$$
\int_Mz(B(s_1, s_2), s_3),
$$
where $z$ is the pairing $z:(T(V_3), \ T_c(V_3^+))\to\Vol$, and $V^+:=V^*\otimes \tr$. The displayed
formula makes sense if the supports of $s_1$, $s_2$ and $s_3$ belong
to one neighborhood $\cU$; due to the locality of the operators this
suffices.

Let the weight of the representation $\rho^*_3$ be
$$
(n_1, n_2)=(l_1+m_1-s, l_2+m_2+s-d); \quad \nu=n_1-n_2=\lambda+\mu-2s+d.
$$
Since the dimensions of the spaces
$V$ and $V^+$ are equal, it suffices to look for operators such that
$\lambda\leq\mu\leq\nu$ (the other
operators will be 1-dual or 2-dual to such operators).

Thus, let us solve our
system under the condition that $\lambda\leq\mu\leq\lambda+\mu-2s+d$, i.e.,
$$
\mu\geq\lambda\geq 2s-d.
$$
             
\begin{Lemma}\label{L4.2} The system
$$
\begin{array}{l}
 a_jx'_-u_j+b_jx''_-u_{j+1}=0, \text{ where}\\
 j=1, 2, \dots, d; a_j\neq0, b_j\neq0,
u_j\in V_1^*\otimes V_2^*
\end{array}
$$
has no solutions of the form $(4)$ if
$s>d$ and $\lambda, \mu\geq 2s-d$.
\end{Lemma}

\begin{proof} Having multiplied each $u_j$ by the corresponding
coefficient we may assume that ${a_j=b_j=1}$. Let $u_j=(s, P_j(i))$,
where $\deg P_j\leq d$. By formulas (5), (6), we have
$$
i(i-\lambda-1)P_j(i-1)=(i+\mu-s)(i-s-1)P_{j+1}(i)
\quad \textup{for all} \  
i \in \{ 0, 1, \dotsc, s+1\}, \ j \in \{1, 2, \dotsc, d\}. 
$$
Observe that $(i+\mu-s)(i-s-1)\neq0$ for $0\leq i\leq d$.

For $i=0$, we get
$$
0=P_{j+1}(0)(\mu-s)(-s-1)
$$
implying that $P_j=0$ for
$j=2, 3, \dots, d$.

For $i=1$, we get
$$
P_j(1)=\frac{1\cdot(1-\lambda-1)}{(1+\mu-s)(1-s-1)}
P_{j-1}(0)=0
$$
for $j=3, 4, \dots, d$.
Substituting $i=2, 3, \dots, d-1$ we get
$$
\text{$P_d(i)=0$ for $i=0, 1, \dots, d-2; \quad
P_{d+1}(i)=0$ for $i=0, 1, \dots, d-1.$}
$$
Moreover, substituting $i=s+1$ in the last equation we get $P_d(s)=0$.
Since $P_d(i)$ and $P_{d+1}(i)$ are of degree $d$, it follows that
$$
P_d(i)=ai(i-1)\dots (i-d+2)(i-s);\quad
P_{d+1}(i)=bi(i-1)\dots (i-d+1).
$$
Having substituted
$$
\begin{array}{l}
 i(i-\lambda-1)a(i-1)(i-2)\dots(i-d+1)(i-s-1)=\\
(i+\mu-s)(i-s-1)bi(i-1)\dots(i-d+1)\end{array}
$$
in the last equation we deduce for $i=d, d+1$ that
$$
a(i-\lambda-1)=b(i+\mu-s).
$$
Therefore,
either $a=b=0$ or $s=\lambda+\mu+1$; the latter contradicts the conditions
$s>d$ and $\lambda, \mu\geq 2s-d$.
\end{proof}

\section{The solutions of degree $d=0$}

All the vectors of degree 0 are annihilated by $\cL_1$ and the
$L_0$-action on them coincides with ${\rho_1^*\otimes\rho_2^*}$.
Therefore, the solutions of system (3) are all the highest weight
vectors from $\cI_0(V_1^*, V_2^*)=1\otimes(V_1^*\otimes V_2^*)$. To
find them, we have to decompose the representation
$\rho_1^*\otimes\rho_2^*$ into the sum of irreducible representations.
This is a~classical problem (for its solution in some cases and an
algorithm, see Table 5 in \cite{OV}). The embedding $V_3^*\to
V_1^*\otimes V_2^*$ generates the map
$$
Z^*: I(V_3^*)\to I(V_1^*)\otimes I(V_2^*)
$$
and the dual projection
$V_1\otimes V_2\to V_3$ gives rise to the
operator
$$
Z: T(V_1)\otimes T(V_2)\to T(V_3).
$$
The above arguments hold for any $n=\dim M$, even for supermanifolds
and any arity of the operator, not only binary.

\section{Solutions of degree $d=1$}

The generic degree 1 element is of the form
$$
f=\partial'_1u_1+\partial'_2u_2+\partial''_1u_3+\partial''_2u_4.
$$
We have
$$
X_+f=-\partial'_2u_1+\partial'_1(X_+u_1)+\partial'_2(X_+u_2)-
\partial''_2u_3+\partial''_1(X_+u_3)+\partial''_2(X_+u_4)
$$
wherefrom
$$
u_1=X_+u_2, \quad X_+u_1=0, \quad u_3=X_+u_4, \quad X_+u_3=0
$$
or
$$
(X_+)^2u_2=(X_+)^2u_4=0.
$$
Hence, $u_2$ and $u_4$ are of the form (4). The remaining two equations yield
$$
(x_2^2\partial_2)f=2h'_2u_2+2h''_2u_4=0, \quad
(x_2^2\partial_1)f=2x'_-u_2+2x''_-u_4=0
$$
wherefrom Lemma \ref{L4.2} implies that for $s\geq 2$ there are no solutions such that
$\lambda\leq\mu\leq\nu$. It remains to consider the cases $s=0, 1$.

Let $s=1$. Then, $\mu\geq\lambda\geq 2s-d=1$. The generic form of the elements
$u_2$ and $u_4$ is
$$
u_2=av_0\otimes w_1+bv_1\otimes w_0 \text{ shortly } a(01)+b(10)
$$
and
$$
u_4=\alpha(01)+\beta(10).
$$
We have to
find all the $u_2$ and $u_4$ satisfying
$$
x'_-u_2+x''_-u_4=0, \quad h'_2u_2+h''_2u_4=0.
$$
Consider the following 3 cases:

1) $\lambda\geq 2$, $\mu\geq 2$. Then, 
$$
x'_-u_2+x''_-u_4=a(11)+2b(20)+2\alpha(02)+\beta(11)=0
$$ implying
$\alpha=b=0$, $\beta=-a$. We have
$$
h'_2u_2+h''_2u_4=la(01)-ma(10)
$$
implying
$l=m=0$.
Thus,
$$
\bar\lambda=(\lambda, 0), \quad \bar\mu=(\mu, 0), \quad
\bar\nu=(\lambda+\mu-s, 0+0+s-1)=
(\lambda+\mu-1, 0).
$$
The corresponding operator is the Schouten concomitant (the
operator $P_3$ on our list).

2) $\lambda=1$, $\mu\geq 2$. We have
$$
x'_-u_2+x''_-u_4=
a(11)+2\alpha(02)+\beta(11)\Longrightarrow\alpha =0, \beta=-a,
$$
then
$$
h'_2u_2+h''_2u_4=la(01)+(l+1)b(10)-ma(10)=0
$$
implying
$la=0, (l+1)b-ma=0$.
There are two cases:

a) $a=0$ and then $b\neq 0$ (since otherwise $f=0$)
and $l=-1$. We have
$$
\bar\lambda=(0, -1), \bar\mu(m_1, m_2), \bar\nu=(m_1-1, m_2-2).
$$
The corresponding operator is $B(w, s)=dw\circ s$ (the operator $P_1$ on our list).

b) $a\neq 0\Longrightarrow l=0, b=ma$. Hence,
$$
\bar\lambda=(1, 0),
\bar\mu=(m_1, m_2), \bar\nu=(m_1, m_2).
$$
The corresponding operator is
$B(\xi, s)=L_{\xi}s$, the Lie derivative
(the operator $P_2$ on our list).

3) $\lambda=\mu=1$. We have
$$
x'_-u_2+x''_-u_4=a(11)+\beta(11)=0\Longrightarrow\beta =-a,
$$
and we have
$$
h'_2u_2+h''_2u_4=la(01)+(l+1)b(10)+(m+1)\alpha(01)-ma(1, 0)=0
$$
or, equivalently,
$$
\left\{\begin{matrix}l
a+(m+1)\alpha=0\cr
(l+1)b-ma=0\cr
\end{matrix}\right .
$$
implying $a=(l+1)(m+1)x$, $b=m(m+1)x$, $\alpha=-l(l+1)x$.

Therefore,
$$
\bar\lambda=(l+1, l),\quad \bar\mu=(m+1, m)\quad \bar\nu=(l+m+1, l+m)
$$
that is the
corresponding operator is of type $P_4$ in our notations.

The case $m=0$, $l=-1$ as well as $l=0$, $m=-1$ and $l=m=-1$ are particular
cases which correspond to {\it two} operators each:
$$
B(\xi, w)=a\xi dw+bL_{\xi}w
$$
in the first two cases and
$$
B(w_1, w_2)=aw_1dw_2+bw_2dw_1
$$
in the third case. (All these operators correspond
to operators $P_1$ and $P_2$ on our list).

Let now $s=0$, $u_2=a(00)$,
$u_4=b(00)$. We have:

1) $\lambda\geq 1$, $\mu\geq 1$, then
$$
x'_-u_2+x''_-u_4=a(10)+b(01)=0\Longrightarrow a=b=0.
$$

2) $\lambda=0$, $\mu\geq 1$, then
$$
x'_-u_2+x''_-u_4=b(01)=0\Longrightarrow b=0
$$
and
$$
h'_2u_2+h''_2u_4=la(00)\Longrightarrow l=0.
$$
Since
$$
\bar\lambda=(00), \quad \bar\mu(m_1, m_2), \quad \bar\nu(m_1, m_2-1),
$$
the corresponding
operator is $P_1$.

3) $\lambda =\mu=0$. We see that the condition $x'_-u_2+x''_-u_4=0$ holds always.
The other condition takes the
form
$$
h'_2u_2+h''_2u_4=la(00)+mb(00).
$$

If $l$ and $m$ do not vanish simultaneously, then $a=mx$ and
$b=-lx$. We have
$$
\bar\lambda=(l, l), \quad \bar\mu=(m, m), \quad \bar\nu=(l+m, l+m-1),
$$
hence,
the corresponding operator is of type $P_4$.

If $l=m=0$, then we have two operators of type $P_1$:
$$
B(f, g)=afdg+bgdf.
$$
Thus, we have verified that for $n=2$ all the first
order operators are listed.

\section{Solutions of degree 2}

The generic form of a~degree 2 vector is
$$

\begin{array}{l}
f=\partial'{}^2_1u_1+\partial'_1\partial'_2u_2+\partial'{}^2_2u_{3}+
\partial'_1\partial''_1u_4+\partial'_1\partial''_2u_5+\\
\partial'_2\partial''_1u_6+
\partial'_2\partial''_2u_7+\partial''_2\partial'_1u_8+\partial''_1\partial''_2u_9+
\partial''{}^2_2u_{10}.
\end{array}
$$
We have
$$

\begin{array}{l}
X_+f=-2\partial'_1\partial'_2u_1+\partial'{}^2_1(X_+u_1)-\partial'{}^2_2u_2+
\partial'_1\partial'_2(X_+u_2)+\partial'{}^2_2(X_+u_3)-\partial'_2\partial''_1u_4-\\
\partial'_1\partial''_2u_4+\partial'_1\partial''_1(X_+u_4)-
\partial'_2\partial''_2u_5+\partial'_1\partial''_2(X_+u_5)-
\partial'_2\partial''_2u_6+\partial'_2\partial''_1(X_+u_6)+\\
\partial'_2\partial''_2(X_+u_7)-2\partial''_1\partial''_2u_8+
\partial''{}^2_1(X_+u_8)-\partial''{}^2_2u_9+\partial''_1\partial''_2(X_+u_9)+
\partial''{}^2_2(X_+u_{10})=0
\end{array}
$$
which implies
$$
\begin{matrix}X_+u_1=0, & X_+u_4=0, & X_+u_8=0, \cr
X_+u_2-2u_1=0, & X_+u_5-u_4=0, & X_+u_9-2u_8=0, \cr
X_+u_3-u_2=0, & X_+u_6-u_4=0, & X_+u_{10}-u_9=0, \cr
& X_+u_7-u_5-u_6=0.& \cr
\end{matrix}
$$
All
these vectors can be expressed in terms of $u_3$, $u_7$, $u_{10}$ and $u_0$:
$$
\begin{matrix}
u_1=\frac{1}{2}(X_+)^2u_3,& u_4=\frac{1}{2}(X_+)^2u_7,&
u_8=\frac{1}{2}(X_+)^2u_{10},\cr
u_2=X_+u_3,& u_5=\frac{1}{2}X_+u_7-u_0,& u_9=X_+u_{10},\cr
 &u_6=\frac{1}{2}X_+u_7+u_0.& \cr
\end{matrix}
$$
Moreover,
$$
(X_+)^3u_3=(X_+)^3u_7=(X_+)^3u_{10}=0, \qquad X_+u_0=0.
$$
The condition
$(x_2^2\partial_1)f=0$ implies
$$

\begin{array}{l}
2\partial'_1(x'_-u_2)-2\partial'_1u_3+4\partial'_2(x'_-u_3)+
2\partial'_1(x''_-u_5)+2\partial''_2(x'_-u_6)+2\partial''_2(x''_-u_7)+\\
2\partial''_1(x''_1u_9)-2\partial''_1u_{10}+4\partial''_2(x''_-u_{10})=0
\end{array}
$$
wherefrom
$$
2x'_-u_3+x''_-u_7=0\eqno {(7)}
$$
$$
x'_-u_7+2x''_-u_{10}=0 \eqno {(8)}
$$
$$
x'_-u_2+x''_-u_5-u_3=0 \eqno {(9)}
$$
$$
x'_-u_6+x''_-u_9-u_{10}=0 \eqno{(10)}
$$
Equations (7), (8) imply thanks to Lemma \ref{L4.2} that for $s\geq 3$
and $\mu\geq\lambda\geq 2s-2$
$$
u_3=u_7=u_{10}=0.
$$
But then formula (9) implies that $x''_-u_0=0$ and (10)
implies $x'_-u_0=0$. Hence,
$$
u_0=cv_{\lambda}\otimes w_{\mu}
$$
but
$X_+(v_{\lambda}\otimes w_{\mu})\neq 0$
for $\lambda$ and $\mu$ indicated; hence, $u_0=0$, i.e., $f=0$.

There remain the cases $s=0$, 1, 2.

Let $s=2$. Then, $\mu\geq\lambda\geq 2s-d=2$. The equation
$x_2^2\partial_2f=0$ implies
$$
h'_2u_2+h''_2\cdot u_5=0\eqno {(11)}
$$
$$
(2h'_2-1)u_3+h''_2u_7=0 \eqno{(12)}
$$
$$
h'_2u_6+h''_2u_9=0 \eqno {(13)}
$$
$$
h'_2u_7+(2h''_2-1)u_{10}=0 \eqno{(14)}
$$
The generic form of the elements $u_3, u_7, u_{10}, u_0$ is as
follows:
$$
\begin{matrix}
u_3=a\cdot(02)+b(11)+c(20),& u_{10}=x(02)+y(11)+z(20),\cr
u_7=2(\alpha(02)+\beta(11)+\gamma(20),& u_0=p(01)+q(10).\cr
\end{matrix}
$$
Therefore,
$$
\begin{matrix}
 u_2&=((\mu-1)a+\lambda b)(01)+(\mu b+(\lambda -1)c)(10),\cr
u_5&=((\mu -1)\alpha +\lambda\beta
-p)(01)+(\mu\beta+(\lambda-1)\gamma-q)(10),\cr
u_6&=((\mu -1)\alpha +\lambda\beta +p)(01)+(\mu\beta+(\lambda-1)\gamma+q)(10),\cr
u_9&=((\mu-1)x+\lambda y)(01)+(\mu y+(\lambda-1)z)(10).
\end{matrix}
$$
Let us substitute $u_j$ in the system (7)--(8). We get
$$
-a+2(\mu-1)\alpha+2\lambda\beta-2p=0\quad \text{for $\mu\geq 2$}\eqno{(15)}
$$
$$
(\mu-1)a+(\lambda-1)b+\mu\beta+(\lambda-1)\gamma -q=0 \quad\text{
for $\lambda\geq 1$, $\mu\geq 1$}\eqno {(16)}
$$
$$
2\mu b+(2\lambda -3)c=0\quad \text{ for $\lambda\geq 2$} \eqno{(17)}
$$
$$
(2\mu -3)x+2\lambda y=0 \quad \text{ for $\mu\geq 2$} \eqno {(18)}
$$
$$
(\mu -1)\alpha +\lambda\beta+p +(\mu-1)y+(\lambda-1)z=0
\quad \text{ for $\lambda\geq 1$, $\mu\geq 1$} \eqno {(19)}
$$
$$
2\mu\beta+2(\lambda-1)\gamma+2q-z=0 \quad \text{
for $\lambda\geq 2 $}\eqno {(20)}
$$
One more equation is obtained from the condition $X_+u_0=0$, namely,
$$
\mu p+\lambda q=0.\eqno {(21)}
$$
Substituting $u_j$ into the system (9)--(10) we get
$$
\begin{array}{ll}
 \alpha=0, \quad x=0& \text{ for }\mu\geq 3\cr
a+2\beta=0, \quad \alpha +2y=0& \text{ for }\lambda\geq 1, \; \mu\geq 2\cr
2b+\gamma=0, \quad 2\beta +z=0& \text{ for }\lambda\geq 2, \; \mu\geq 1\cr
c=0, \quad \gamma =0& \text{ for }\lambda \geq 3\cr
\end{array}\eqno {(22)}
$$
Consider the 3 cases:

1) $\lambda\geq 3$, $\mu\geq 3$. Then, (22) implies that
$b=c=\alpha=\gamma=x=y=0$ and $ a=z=-2\beta$ and we have
$$
\begin{matrix}
 &\left.\begin{matrix}
(15)&\Longrightarrow& 2\beta +2\lambda\beta -2p=0\cr
(19)&\Longrightarrow&\lambda\beta -2(\lambda -1)\beta +p=0\cr
\end{matrix}
\right\}&\Longrightarrow\beta =p=0,\cr
(21)&\Longrightarrow& \mu p+\lambda q=0&\Longrightarrow q=0.
\end{matrix}
$$
No solutions.

2) $\lambda =2$, $\mu\geq 3$. Then, (22) implies
$\alpha =x=y=0$, $a=z=-2\beta$, $\gamma =-2b$ and we have
$$
\begin{matrix}
 (19)\Longrightarrow&+2\beta +p -2\beta =0& \Longrightarrow p=0,\cr
(21)\Longrightarrow&\mu p+2q=0&\Longrightarrow q=0,\cr
(15)\Longrightarrow& 2\beta +4\beta -2p=0&\Longrightarrow\beta =0,\cr
(16)\Longrightarrow& -2(\mu -1)\beta +b+\mu\beta -2b-q=0&\Longrightarrow b=0,\cr
(17)\Longrightarrow& 2\mu b+c=0&\Longrightarrow c=0.\cr
\end{matrix}
$$
No solutions again.

3) $\lambda =\mu =2$. Then, (22) implies $a=z=-2\beta$, $\gamma=-2b$,
$\alpha =-2y$, and we have
$$
\begin{matrix}
 (15) \Longrightarrow&2\beta -4y+4\beta -2p=0,\cr
(16) \Longrightarrow&-2\beta +b+2\beta -2b-q=0,\cr
(17) \Longrightarrow&4b+c=0,\cr
(18)\Longrightarrow&x+4y=0,\cr
(19) \Longrightarrow&-2y+2\beta +p+y-2\beta =0,\cr
(20) \Longrightarrow&4\beta -4b +2q +2\beta =0,\cr
(21) \Longrightarrow&2p +2q=0.
\end{matrix}
$$
The solution of this system is
$$
\beta =b=p=y=-q, \qquad c=x=4q, \qquad a=\alpha =\gamma =z=2q.
$$
But having substituted
$u_3=2q(02)-q(11)+4q(20)$, and $u_7=4q(02)-2q(11)+4q(20)$ into (12) we get
$$

\begin{array}{l}
2(2l-1)q(02)-(2l+1)q(11)+4(2l+3)q(20)+\\
4(m+2)q(02)+2(m+1)q(11)+4mq(20)=0
\end{array}
$$
implying $q=0$.

Now let $s=1$. Then, 
$$
\begin{matrix}
u_3=a(01)+b(10), & u_7=2(\alpha (01)+\beta (10))&
u_{10}=x(01)+y(10)\cr
u_2=(\mu a+\lambda b)(00)& u_5=(\mu\alpha +\lambda\beta -p)(00)&
u_9=(\mu x+\lambda y)(00)\cr & u_6=(\mu\alpha +\lambda\beta +p)(00)& \cr
\end{matrix}
$$
Let us substitute this into (7)--(8). We get
$$
-a+\mu\alpha +\lambda\beta -p=0 \quad\text{for }\mu\geq 1\eqno{(23)}
$$
$$
\mu a+(\lambda -1)b=0 \quad\text{for }\lambda\geq 1\eqno{(24)}
$$
$$
(\mu -1)x+\lambda y=0 \quad\text{for }\mu\geq 1\eqno{(25)}
$$
$$
\mu\alpha +\lambda\beta +p-y=0 \quad\text{for }\lambda\geq 1\eqno{(26)}
$$
The system
(9)--(10) yields
$$
\begin{array}{lll}
\alpha =0, &x=0&\text{for }\mu\geq 2\cr
a+\beta=0, &x+y=0&\text{for }\lambda\geq 1, \; \mu\geq 1\cr
b=0, &\beta =0 &\text{for }\lambda\geq 2\cr
\end{array}\eqno {(27)}
$$

1) Let $\lambda\geq 2$, $\mu\geq 2$. Then, equation (27) implies
$a=b=\alpha =\beta =x=y=0$ and (23) implies $p=0$. No solutions.

2) $\lambda =1$, $\mu\geq 2$. From (27) it follows that
$\alpha =x=y=0$, $\beta =-a$ and equations (24) and (23) imply that
$(24)\Longrightarrow a=0$ and $(23)\Longrightarrow p=0$,
respectively.

Having substituted $u_3=b(10)$, $u_2=b(00)$, and $u_5=u_7=0$ into (11) and (12)
we get $(l+1)b=(2l-1)b=0$. Hence, $b=0$. No solutions.

3) $\lambda =0$, $\mu\geq 2$. Equation (27) implies $\alpha =x=0$. Moreover,
$b=\beta =y=0$ (since there is no vector $v_1$). From
(23) we get $p=a$. Having substituted $u_3=a(01)$, $u_7=0$, $u_6=a(00)$, $u_9=0$
into (12), (13) we get $(2l-1)a=0,$ $la=0\Longrightarrow a=0$. No solutions.

4) $\lambda =\mu =1$. Equation (27) implies $\beta =-a$, $\alpha =-y$.
Having substituted this into (23)--(26) we get
$$
\begin{matrix}
 (23)\Longrightarrow& -a-y-a-p=0, \cr
(24)\Longrightarrow &a=0, \text { hence } a=y=p=0, \cr
(25)\Longrightarrow &y=0, \cr
(26)\Longrightarrow &-y-a+p-y=0.
\end{matrix}$$
Having substituted $u_3=b(10)$, $u_2=b(00)$, $u_{10}=x(01)$, $u_9=x(00)$,
$u_5=u_6=u_7=0$ into (11)--(14) we get
$$
\left.\begin{matrix}
(11)& lb(00)=0\cr
(12)& (2l+1)b(10)=0\cr
(13)& mx(00)=0\cr
(14)& (2m+1)x(01)=0\cr
\end{matrix}\right\}\Longrightarrow b=x=0.
$$

5) $\lambda =0$, $\mu=1$. Then, $b=\beta y=0$ (since there is no vector $v_1$).
The system (27) does not give
anything new. From (23)--(26) we deduce that $p=\alpha -a$.
Thus, $u_3=a(01)$, $u_2=a(00)$, $u_7=2\alpha (01)$, $u_5=a(00)$,
$u_6=(2\alpha-a)(00)$, $u_{10}x(01)$, $u_9=x(00)$. Having substituted this
into (11)--(14) we get
$$
\begin{matrix}
(11) & la(00)+ma(00)=0,\cr
(12)& (2l-1)a(01) +2(m+1)\alpha (01)=0, \cr
(13)& l(2\alpha -a)(00) + mx(00)=0, \cr
(14) &2l\alpha (01)+(2m+1)(01)=0.
\end{matrix}
$$

a) Let $a\neq 0$. Then, from the first equation, i.e., (11), we get $-l=m$; hence,
$$
\begin{matrix}
 (2l-1)a+2(-l+1)\alpha =0, \cr
l(2\alpha -a-x)=0, \cr
2l\alpha +(-2l+1)x=0.
\end{matrix}
$$
From the first equation
we get $a=(l-1)c$, $\alpha =(l-\frac{1}{2})c$. From the second equation (or
from the third one if $l=0$) we get $x=lc$.

Thus,
$$
\bar\lambda =(l, l),\quad \bar\mu =(-l+1, l),\quad \bar\nu =(0, -1),
$$
i.e., the corresponding operator is of the form $S^{*1}_2$

b) $a=0$. Then, 
$$
(m+1)\alpha =0, \quad
2l\alpha +mx =0, \quad
2l\alpha +(2m+1)x=0.
$$

If $\alpha=0$, then from the second and the third equations we derive
$x=0$, i.e., $f=0$. Hence, $\alpha\neq 0$, $m=-1$, $x=2l\alpha$.

Thus,
$$
\bar\lambda =(l, l),\quad \bar\mu =(0, -1),\quad \bar\nu =(l-1, l-2)
$$
and the
corresponding operator is of type $S_2$.

6) $\lambda =\mu =0$. Then, $u_3=u_7=u_{10}=0$ but $u_0\neq 0$. The equations
(23)--(27) do not give anything. Let us substitute $u_5=-p(00)$, $u_6=p(00)$
into (11)--(13):
$$
\left.\begin{array}{ll}(11) &-mp(00)=0\cr
(13) &lp(00)=0\cr
\end{array}\right\}\Longrightarrow l=m=0.
$$

Thus, $\bar\lambda =\bar\mu =(0, 0)$, $\bar\nu =(-1, -1)$,
and the corresponding operator $B(f, g)=df\wedge dg$ is of
type $S_1$.

There remains the case $s=0$. In this case,
$u_3=a(00)$, $u_7=2\alpha (00)$, $u_{10}=x(00)$ all the other $u_j$ being zero.
From (9), (10) we deduce that $u_3=u_{10}=0$. Hence, there
remains only $u_7=2\alpha (00)$. From (12), (14) we see that
$m\alpha =k\alpha =0$, and therefore 
$$
\bar\lambda =\bar\mu =(0, 0),\quad \bar\nu =(0, -2).
$$

The corresponding
operator $B(f, g)=Z(df; dg)$ is of type $S_1$.

\section{The solutions of degree 3}

The generic form of a~homogeneous element of degree 3 is
$$

\begin{array}{l}
f=\partial'{}^3_1u_1+\partial'{}^2_1\partial'_2u_2+\partial'_1\partial'{}^2_2u_3+
\partial'{}^3_2u_4+\partial'{}^2_1\partial''_1u_5+\\
\partial'_1\partial'_2\partial''_1u_6+\partial'{}^2_2\partial''_1u_7+
\partial'{}^2_1\partial''_2u_8+\partial'_1\partial'_2\partial''_2u_9+
\partial'{}^2_2\partial''_2u_{10}+\partial'_1\partial''{}^2_1u_{11}+\\
\partial'_1\partial''_1\partial''_2u_{12}+\partial'_1\partial''{}^2_2u_{13}+
\partial'_2\partial''{}^2_1u_{14}+\partial'_2\partial''_1\partial''_2u_{15}+
+\partial'_2\partial''{}^2_2u_{16}+\partial''{}^3_1u_{17}+\\
\partial''{}^2_1\partial''_2u_{18}+
\partial''_1\partial''{}^2_2u_{19}+\partial''{}^3_2u_{20}.
\end{array}
$$
The equation $X_+f=0$
yields
$$

\begin{matrix}
X_+u_1=0& X_+u_5=0& X_+u_{11}=0& X_+u_{17}=0\cr
X_+u_2=3u_1& X_+u_6=2u_5& X_+u_{12}=2u_{11}& X_+u_{18}=3u_{17}\cr
X_+u_3=2u_2& X_+u_7=u_6& X_+u_{13}=u_{12}& X_+u_{19}=2u_{18}\cr
X_+u_4=u_3& X_+u_8=u_5& X_+u_{14}=u_{13}& X_+u_{20}=u_{19}\cr
 &X_+u_9=u_6+2u_8& X_+u_{15}=u_{12}+2u_{14}& \cr
 &X_+u_{10}=u_7+u_9& X_+u_{16}=u_{13}+u_{15}& \cr
 \end{matrix}
 $$
 The solution of this system is
$$

\begin{matrix}
u_1=\frac{1}{6}(X_+)^3u_4& u_5=\frac{1}{6}(X_+)^3u_{10}&
u_{11}=\frac{1}{6}(X_+)^3u_{16}& u_{17}=\frac{1}{6}(X_+)^2u_{20}\cr
u_2=\frac{1}{2}(X_+)^2u_4& u_6=\frac{1}{3}(X_+)^2u_{10}-X_+u'&
u_{12}=\frac{1}{3}(X_+)^2u_{16}-X_-u''& u_{18}=\frac{1}{2}(X_+)^2u_{20}\cr
u_3=X_+u_4& u_7=\frac{1}{3}X_+u_{10}-u'& u_{13}=\frac{1}{3}X_+u_{16}-u''&
u_{19}=X_+u_4\cr
 & u_8=\frac{1}{6}(X_+)^2u_{10}+X_+u'& u_{14}=\frac{1}{6}(X_+)^2u_{16}+X_+u''\cr
 & u_9=\frac{2}{3}X_+u_{10}+u'& u_{15}=\frac{2}{3}X_+u_{16}+u''\cr
 \end{matrix}
 $$
where
$$
(X_+)^4u_4=(X_+)^4u_{10}=(X_+)^4u_{16}=(X_+)^4u_{20}=0,
$$
$$
(X_+)^2u'=
(X_+)^2u''=0.
$$
From the equation $(x^2_2\partial_1)f=0$ we get
$$
x'_-u_2-u_3+x''_-u_8=0 \eqno {(28)}
$$
$$
2x'_-u_3-3u_4+x''_-u_9=0 \eqno {(29)}
$$
$$
3x'_-u_4+x''_-u_{10}=0 \eqno {(30)}
$$
$$
x'_-u_6-u_7+x''_-u_{12}-u_{13}=0 \eqno {(31)}
$$
$$
2x'_-u_7+x''_-u_{15}-u_{16}=0 \eqno {(32)}
$$
$$
x'_-u_9-u_{10}+2x''_-u_{13}=0 \eqno {(33)}
$$
$$
x'_-u_{10}+x''_-u_{16}=0 \eqno {(34)}
$$
$$
x'_-u_{14}+x''_-u_{18}-u_{19}=0 \eqno {(35)}
$$
$$
x'_-u_{15}+2x''_-u_{19}-3u_{20}=0 \eqno {(36)}
$$
$$
x'_-u_{16}+3x''_-u_{20}=0 \eqno {(37)}
$$
From equations (30), (34), (37) thanks to Lemma \ref{L4.2} it follows that
\[
    u_4=u_{10}=u_{16}=u_{20}=0
    \quad \textup{for }\ s\geq 4, \ \lambda\leq\mu\leq\nu.
\]
Thanks to the same lemma equation (32) implies $u'=u''=0$.

From the equation
$(x_2^2\partial_2)f=0$ we get
$$
h'_2u_2+h''_2u_8=0 \eqno {(38)}
$$
$$
(2h'_2-1)u_3+h''_2u_9=0 \eqno {(39)}
$$
$$
(3h'_2-3)u_4+h''_2u_{10}=0 \eqno {(40)}
$$
$$
h'_2u_6+h''_2u_{12}=0 \eqno {(41)}
$$
$$
(2h'_2-1)u_7+h''_2u_{13}=0 \eqno {(42)}
$$
$$
h'_2u_9+(2h''_2-1)u_{15}=0 \eqno {(43)}
$$
$$
(2h'_2-1)u_{10}+(2h''_2-1)u_{16}=0 \eqno {(44)}
$$
$$
h'_2u_14 +h''_2u_{18}=0 \eqno {(45)}
$$
$$
h'_2u_{15}+(2h''_2-1)u_{19}=0 \eqno {(46)}
$$
$$
h'_2u_{16}+(3h''_2-3)u_{20}=0 \eqno {(47)}
$$
Let $s=3$. Denote
$$
u_{j_i}=a_i(03)+b_i(12)+c_i(21)+d_i(30)\text{ for
$(j_1=4, j_2=10, j_3=16, j_4=20)$}.
$$
From (30), (34), (37) we get
$$
\begin{matrix}
a_2=a_3=a_4=0\text{for }\mu\geq 4\cr
3a_1+3b_2=a_2+3b_3=a_3+9b_4=0\text{for }\lambda\geq 1, \; \mu\geq 3\cr
6b_1+2c_2=2b_2+2c_3=2b_3+6c_4=0\text{for }\lambda\geq 2, \; \mu\geq 2\cr
9c_1+d_2=3c_2+d_3=3c_3+3d_4=0\text{for }\lambda\geq 3, \; \mu\geq 1\cr
d_1=d_2=d_3=0\text{for }\lambda\geq 4\cr
\end{matrix}\eqno {(48)}
$$
One more system
is obtained from (28), (29), where
$$
\begin{matrix}
u'=p(02)+q(11)+r(20), \\
u_4=a(03)+b(12)+c(21)+d(30);\\
u_{10}=3(\alpha(03)+\beta(12)+\gamma(21)+\delta(30))
-(\mu-2)a-\lambda b+\\
(\mu-1)(\mu-2)\alpha+2\lambda(\mu-1)\beta+
\lambda(\lambda-1)\gamma+2(\mu-1)p+2\lambda q=0 \; \text{ for }\mu\geq 2
\end{matrix}\eqno {(49)}
$$
$$
\begin{matrix}
\frac{1}{2}(\mu-1)(\mu-2)a+(\lambda-1)(\mu-1)b+\frac{1}{2}(\lambda-1)
(\lambda-2)c+\\
\frac{1}{2}\mu(\mu-1)\beta +(\lambda-1)\mu\gamma +\frac{1}{2}
(\lambda -1)(\lambda -2)\delta +\mu q+(\lambda -1)r=0\text{ for }\lambda\geq 1,
\mu\geq 1 \end{matrix}\eqno {(50)}
$$
$$
\mu (\mu -1)b+(2\lambda -3)\mu c+(\lambda -2)^2d=0 \text{ for }\lambda\geq 2 \eqno
{(51)}
$$
$$
-3a+6(\mu -2)\alpha +6\lambda\beta +3p=0\text{ for }\mu\geq 3 \eqno {(52)}
$$
$$
2(\mu -2)a+(2\lambda -3)b+4(\mu -1)\beta +4(\lambda -1)\gamma +2q=0
\text{ for }\mu\geq 2, \lambda\geq 1 \eqno {(53)}
$$
$$
4(\mu -1)b+(4\lambda -7)c+2\mu\gamma +2(\lambda -1)\delta +r=0
\text{ for }\mu\geq 1, \lambda\geq 2 \eqno {(54)}
$$
$$
6\mu c+(6\lambda -15)d=0\text{ for }\lambda\geq 3 \eqno {(55)}
$$
And one more equation from the
condition $(X_-)^2u'=0$:
$$
\mu (\mu -1)p+2\lambda\mu q+\lambda (\lambda -1)r=0. \eqno {(56)}
$$
Let us consider the corresponding cases.

1) $\lambda\geq 4$, $\mu\geq 4$. From (48) we get
$b_1=c_1=d_1=a_2=c_2=d_2=0$, $b_2=-a_1$, $u_4=a(03)$, $u_{10}=-a(12)$.
Having substituted this into the system (49)--(56) we
get: $a=-3\beta$, $b=c=d=\alpha =\gamma =\delta =0$ and
$$
\begin{matrix}
(52) &9\beta +6\lambda\beta +3p=0\cr
(53) &-6(\mu -2)\beta +4(\mu -1)\beta +2q=0\cr
(54) &r=0\cr
(56) &\mu (\mu -1)p+2\lambda\mu q+\lambda (\lambda -1)r=0.
\end{matrix}
$$
The determinant of the system is equal to $6\mu (6\lambda +3\mu -3)$.
This determinant is $\neq 0$
for $\lambda\geq 4$, $\mu\geq 4$. Hence, $u_4=u_{10}=u'=0$.
Similarly, $u_{16}=u_{20}=u''=0$. No solutions.

2) $\lambda =3$, $\mu\geq 4$. Let us first consider the opposite
case: $\lambda\geq 4$, $\mu =3$. From (48) we deduce that
$b_1=c_1=d_1=c_2=d_2=d_3=0$, $a_1=-b_3=-3\beta$. Having substituted
this into (49)--(56) we get
$$
\begin{matrix}
(54) &r=0\cr
(50)&-3\beta +3\beta +3q +(\lambda -1)r=0\Longrightarrow q=0\cr
(56)&6p+6\lambda q +\lambda (\lambda -1)r=0\Longrightarrow p=0\cr
(53)&-6\beta +8\beta +2q =0\Longrightarrow \beta =0\cr
(52)&9\beta +6\alpha +6\lambda\beta +3p=0\Longrightarrow\alpha =0.
\end{matrix}
$$
Thus, $u_4=u_{10}=u'=0$.

Let us return to the case $\lambda =3$, $\mu\geq 4$. By the just proved
$u_{16}=u_{20}=u''=0$. From (48) we derive that $u_4=c(21)+d(30)$,
$u_{10}=-9c(30)$.
Having substituted this into (49)--(56) we get
$$
\begin{matrix}
(52)&p=0\cr
(53) & q=0\cr
(56) &\mu (\mu -1)p+6\mu q+6r=0\Longrightarrow r=0\cr
(54) &5c-12c+r=0\Longrightarrow c=0\cr
(55) &6\mu c+3d=0\Longrightarrow d=0
\end{matrix}
$$
No nonzero solutions.

3) $\lambda =\mu =3$. From system (48) we get
$u_4=a(03)+b(12)+c(21)+d(30)$ and $u_{10}=3(\alpha (03)-\frac{1}{3}a(12)-b(21)
-3c(30))$. Having substituted into (49)--(56) we get
$$
\begin{matrix}
(49) &-a-3b+2\alpha -4a-6b+4p+6q=0\cr
(50)& a+4b+c-a-6b-3c+3q+2r=0\cr
(51) &6b+9c+d=0\cr
(52) &-3a+6\alpha -6a+3p=0\cr
(53) &2a+3b-\frac{8}{3}a-8b+2q=0\cr
(54) &8b+5c-6b-12c+r=0\cr
(55) &18c+3d=0\cr
(56) &6p+18q+6r=0
\end{matrix}
$$
The
system is nondegenerate. No nonzero solutions.

The case $s=3$ is exhausted since $\mu\geq\lambda\geq 2s-d=3$.

Let now $s=2$. In this case $\mu\geq\lambda\geq 1$. For the same $j_{i}$ as for $s=3$, define $u_{j_i}=a_i(02)+b_i(11)+c_i(20)$. Let us substitute
this into equations (30), (34), (37). We get
$$

\begin{array}{rl}
a_2=a_3=a_4=0& \text{ for }\mu\geq 3\\
3a_1+2b_2=a_2+2b_3=a_3+6b_4=0&\text{ for }\lambda\geq 1, \; \mu\geq 2\\
6b_1+c_2=2b_2+c_3=2b_3+3c_4=0 & \text{ for }\lambda\geq 2, \; \mu\geq 1\\
c_1=c_2=c_3=0& \text{ for }\lambda\geq 3.
\end{array}
$$
Let us substitute
$$
u_4=a(02)+b(11)+c(20), \quad u_{10}=3(\alpha (02)+\beta (11)+\gamma (20)),
\quad u'=p(01)+q(10)
$$
into
equations (28), (29). We get
$$
\frac{1}{2}\mu (\mu -1)a+(\lambda -1)\mu b+\frac{1}{2}(\lambda -1)
(\lambda -2)c=0 \text{ for }\lambda\geq 1 \eqno {(58)}
$$
$$
-(\mu -1)a-\lambda b+\frac{1}{2}\mu (\mu -1)\alpha +\lambda\mu\beta +
\frac{1}{2}\lambda (\lambda -1)\gamma +\mu p+
\lambda q=0\text{ for }\mu\geq 1
\eqno {(59)}
$$
$$
4\mu b+(4\lambda -7)c=0 \text{ for }\lambda\geq 2 \eqno {(60)}
$$
$$
2(\mu -1)a+(2\lambda -3)b+2\mu\beta +2(\lambda -1)\gamma +q=0 \text{ for }
\lambda\geq 1, \; \mu\geq 1 \eqno {(61)}
$$
$$
-3a+4(\mu -1)\alpha +4\lambda\beta +2p=0 \text{ for }\mu\geq 2.\eqno {(62)}
$$
Let us consider the corresponding cases.

1) $\lambda\geq 3$, $\mu\geq 3$. Then, (57) implies $u_4=u_{10}=u_{16}=u_{20}=0$.
Therefore,
$$
\left.\begin{matrix}
(61)\Longrightarrow q=0\cr
(62)\Longrightarrow p=0\cr
 \end{matrix}\right \}
\Longrightarrow u'=0.
$$
Similarly, $u''=0$ no solutions.

2) $\lambda =2$, $\mu\geq 3$. Again, let first $\lambda\geq 3$, $\mu
=2$. From (57) we deduce that $u_4=0$, $u_{10}=3\alpha (02)$. Having
substituted this into (58)--(62) we get
$$
\left.\begin{matrix}
(59)& \alpha +2p+\lambda q=0\cr
(61) &q=0\cr
(62) &4\alpha +2p=0 \cr
\end{matrix}
\right \}\Longrightarrow u_4=u_{10}=u'=0.
$$
Let us return to the case $\lambda =2$, $\mu\geq 3$. In this case
$u_{16}=u_{20}=u''=0$ and (57) implies that $u_4=b(11)+c(20)$,
$u_{10}=-6b(20)$. Having substituted this into (58)--(62) we get
$$
\begin{matrix}
(58) &\mu b=0\cr
(59) &-2b-2b +\mu p+2q=0\cr
(61) &4\mu b+c=0\cr
(62) &p=0.
\end{matrix}
$$
No solutions.

3) $\lambda =\mu =2$. From (57) it follows that
$$
u_4=a(02)+b(11)+c(20),\qquad u_{10}=3(\alpha (02)-\frac{1}{2}a(11)-2b(20)).
$$
Having substituted this into the
corresponding equations we get
$$
\begin{matrix}
(58)& a+2b=0\cr
(59)& -a-2b+\alpha -2a-2b+2p+2q=0\cr
(60)& 8b+c=0\cr
(61)& 2a+b-2a-4b+q=0\cr
(62) &-3a+4\alpha -4a +2p=0
\end{matrix}
$$
The solution is
$$
\begin{matrix}
a=\alpha =2x, \quad b=-x, \quad c=8x, \quad p=3x, \quad
q=-3x;\cr
u_4=x(2(02)-(11)+8(20)), \quad u_{10}=x(6(02)-3(11)+6(20)), \quad
u'=3x((01)-(10)).
\end{matrix}
$$
Similarly,
$$
u_{20}=y(8(02)-(11)+2(02)),\quad u_{16}=3y(2(02)-(11)+2(20)),\quad
u''=3y(-(01)+(10)).
$$
Having
substituted $u_{10}$, $u_9=\frac{2}{3}X_+u_{10}+u'=u'$, $u_{13}=
\frac{1}{3}X_+u_{16}-u''=-u''$ into (33) we get
$$
3x(11)-6x(20)-6x(02)+3x(11)-6x(20)+12y(02)-6y(11)=0
$$
implying $x=y=0$.

4) $\lambda =1$, $\mu\geq 3$. There is no vector $v_2$; hence,
$c_1=c_2=c_3=c_4=0$. From (57) we deduce $u_4=a(02)+b(11)$,
$u_{10}=-\frac{3}{2}a(11)$, $u_{16}=u_{20}=0$. Having substituted this
into (58)--(62) we get
$$
\begin{matrix}
(58) &\frac{1}{2}\mu (\mu -1)a=0\cr
(59) &-(\mu -1)a-b-\frac{1}{2}\mu a+\mu p+q=0\cr
(61) &2(\mu -1)a-b-\mu a+q=0\cr
(62) &-3a-2a+p=0.
\end{matrix}
$$
The solution of this system is
$$
a=p=0, \quad q=b, \quad u_4=q(11), \quad u_{10}=0, \quad u'=q(10).
$$
Having substituted
$u_4$, $ u_{10}$, $ u_3=X_+u_4=q(01)+q(10)$, $ u_9=u'$ into (39), (40) we get
$$
\left.\begin{matrix}
(39) &(2l-1)q(01)+\mu q(2l+1)(10)+mq(10)=0\cr
(40) &3lq(11)=0\cr
\end{matrix}\right \}\Longrightarrow q=0.
$$
Moreover, from (31)--(32) it follows that $x''_-(X_+u'')=u''$, $x''_-(X_+u'')=0$,
i.e., $u''=0$.

5) $\lambda =1$, $\mu =2$. From (57) we get $u_4=a(02)+b(11)$,
$u_{10}=3(\alpha (02)-\frac{1}{2}a(11))$ and (58)--(62) become
$$
\begin{matrix}
(58) &a=0\cr
(59) &-a-b+\alpha -a+2p+q=0\cr
(61) &2a-b-2a+q=0\cr
(62) &-5a+4\alpha +2p =0.
\end{matrix}
$$
The solution is
$$
a=\alpha =p=0, \qquad q=b.
$$
To find $u_{16}$, $u_{20}$, $u''$, let us consider the case

$\lambda =2$, $\mu =1$. In this case we find $u_4$, $u_{10}$, $u'$
with the help of equations (57)--(62).

From (57) we deduce that
$u_4=b(11)+c(20)$, $u_{10}=-6b(20)$ and having substituted this into
(58)--(62) we get
$$
\left.\begin{matrix}
(58) &b=0\cr
(59) &-2b+p+2q=0\cr
(60) &4b+c=0\cr
(61) & b-4b+q=0\cr
 \end{matrix}
\right \} \Longrightarrow b=c=p=q=0.
$$
Thus, for $\lambda =2$, $\mu =1$ we have $u_4=u_{10}=u'=0$;
hence, for $\lambda =1$, $\mu=2$ we have
$u_{16}=u_{20}=u''=0$.

The solution $u_4=b(11)$, $u'=b(10)$ does not satisfy
equations (39), (40) as in the case $\lambda =1$, $\mu\geq 3$.

6) $\lambda =\mu 1$. Then, $u_4=b(11)$, $u_{10}=3\beta (11)$ and we have
$$
\begin{matrix}
(58) &0=0\cr
(59) &-b+\beta +p+q=0\cr
(61) &-b+2\beta +q=0.
\end{matrix}
$$
The solution is $p=\beta$, $q=b-2\beta$, $u_4=b(11)$, $u_{10}=3\beta (11)$,
$u'=\beta (01)+(b-2\beta)(10)$.

Similarly, $u_{20}=a(11)$, $u_{16}=3\alpha (11)$,
$u''=(a-2\alpha )(01)+\alpha (10)$.

Let us substitute $u_4$, $u_{10}$, $u_{16}$, $u_{20}$
as well as
$$

\begin{array}{rl}
u_3&=X_+u_4=b(01)+b(10), \cr
u_2&=\frac{1}{2}X_+u_3=b(00), \cr
u_3&=\frac{2}{3}X_+u_{10}+u'=3\beta (01)+b(10), \cr
u_7&=\frac{1}{3}X_+u_{10}-u'=(3\beta -b)(10), \cr
u_8&=\frac{1}{6}(X_+)^2u_{10}+X_+u'=b(00), \cr
u_6&=X_+u_7=(3\beta -b)(00), \cr
u_{19}&=a(01)+a(10), \cr
u_{18}&=a(00), \cr
u_{15}&=a(01)+3\alpha (00), \cr
u_{14}&=a(00), \cr
u_{13}&=(3\alpha -a)(01), \cr
u_{12}&=(3\alpha -a)(00)
\end{array}$$
into the system (38)--(47). We get
$$

\begin{array}{cl}
(38) &lb(00)+mb(00)=0, \cr
(39) &(2l-1)b(01)+(2l+1)b(10)+3(m+1)\beta (01)+mb(10)=0, \cr
(40) &3lb(11)+3(m+1)\beta (11)=0, \cr
(41) & l(3\beta -b)(00)+m(3\alpha -a)(00), \cr
(42) &(2l+1)(3\beta -b)(10)+(m+1)a(01)+3m\alpha (10)=0, \cr
(43) &3l\beta (01) +(l+1)b(10)+(2m+1)(3\alpha -a)(01)=0, \cr
(44) &3(2l+1)\beta (11) +3(2m+1)\alpha (11)=0, \cr
(45) & la(00)+ma(00)=0, \cr
(46) & la(01)+3(l+1)\alpha (10)+(2m+1)a(01)+(2m-1)a(10)=0, \cr
(47) &4(l+1)\alpha (11)+3ma(11)\neq 0.
\end{array}
$$
Consider the two cases:

a) $l+m\neq 0$. Then, (38) and (45) imply that $a=b=0$.
$$

\begin{array}{l}
 (m+1)\beta =0, \\
l\beta +m\alpha =0, \\
(2l+1)\beta +m\alpha =0, \\
l\beta +(2m+1)\alpha =0, \\
(2l+1)\beta +(2m+1)\alpha =0, \\
(l+1)\alpha =0,
\end{array}
$$
where either $\beta\neq 0$ or $\alpha\neq 0$ (i.e., $f\neq 0$).

Let, for example, $\beta\neq 0$. Then, $m=-1$ and $\alpha =l\beta =(2l+1)\beta$
implying $l=-1$.

Thus, $\bar\lambda =\bar\mu =(0, -1)$ and $\bar\nu =(-2, -3)$ and the operator
found is $T_1$.

b) $m=-l$. Then, 
$$

\begin{array}{rlrl}
 (2l-1)b+3(-l+1)\beta =&0& 3l\beta +(-2l+1)(3\alpha -a)=&0\\
(2l+1)b-lb=&0& (l+1)b=&0\\
lb+(-l+1)\beta=&0& (2l+1)\beta +(-2l+1)\alpha =&0\\
l(3\beta -b-3\alpha +a=&0& la+(-2l+1)a=&0\\
(-l+1)a=&0& 3(l+1)\alpha +(-2l+1)a=&0\\
(2l+1)(3\beta -b)-3l\alpha =&0& (l+1)\alpha -la=&0\\
\end{array}
$$
If $l\neq \pm 1$, then $a=b=0$. This case is already considered in heading a).

There remain cases $l=-1$, $ m=1$ and $l=1$, $ m=-1$. Since these
cases are equivalent, let us consider only the first one.

We get
$$
\begin{matrix}
 -3b+6\beta =0& -3\beta +9\alpha -3a=0\cr
0=0& 0=0\cr
-b+2\beta =0& -\beta +3\alpha =0\cr
3\beta -b-3\alpha +a=0& 2a=0\cr
-2a=0& 3a=0\cr
-3\beta +b+3\alpha =0& a=0.\cr
\end{matrix}
$$
The solution of the system is $a=0$, $b=2\beta$, $\alpha =4\beta $.

Thus, $\bar\lambda =(0, -1)$, $\bar\mu =(2, 1)$, $\bar\nu =(0, -1)$ and the
corresponding operator is $T_1^{*1}$.

Let now $s=1$. We get
$$
\begin{matrix}
 u_4=a(01)+b(10)& u_{10}=3(\alpha (01)+\beta (10))\cr
u_3=(\mu a+\lambda b)(00)&
u_9=(2\mu\alpha +2\lambda\beta +p)(00)\cr
u_2=u_8=0, &u'=p(00).
\end{matrix}\eqno{(*)}
$$

Having substituted $(*)$ into the system (28)--(29) we get
$$
-\mu a-\lambda b=0\eqno {(63)}
$$
$$
-a+2\mu\alpha +2\lambda\beta +p=0 \text{ for }\mu\geq 1 \eqno {(64)}
$$
$$
2\mu a+(2\lambda -3)b=0 \text{ for }\lambda\geq 1 \eqno {(65)}
$$
and having
substituted $u_{j_i}=a_i(01)+b_i(10)$ into the system (30), (34), (37) we get
$$
\begin{array}{ll}
a_2=a_3=a_4=0&\text{for }\mu\geq 2\cr
3a_1+b_2=a_2+b_3=a_3+3b_4=0&\text{for }\lambda\geq 1, \; \mu\geq 1\cr
b_1=b_2=b_3=0&\text{for }\lambda\geq 2.\cr
\end{array}\eqno {(66)}
$$

1) $\lambda\geq 2$, $\mu\geq 2$. From (66) we see that
$u_4=u_{10}=u_{16}=u_{20}=0$
and from (64) we get $u'=0$. Similarly, $u''=0$.

2) $\lambda =1$ $\mu\geq 2$. From (66) we get $u_{16}=u_{20}=0$,
$u_4=a(01)+b(10)$, $u_{10}=-3a(10)$. This gives
$$
\left.\begin{matrix}
(63) &-\mu a-b=0\cr
(64) &-a-2a+p=0\cr
(65) &2\mu a-b=0\cr
 \end{matrix}
\right \}\Longrightarrow a=b=p=0.
$$
Moreover, (31) implies $u_{13}=0$, hence, $u''=0$.

3) $\lambda =\mu =1$. From (66) we get $u_4=a(01)+b(10)$,
$u_{10}=3(\alpha (01)-a(10))$. This gives
$$
\left.\begin{matrix}
(63) &-a-b=0\cr
(64) &-a+2\alpha -2a+p=0\cr
(65) &2a-b=0\cr
 \end{matrix}
\right \} a=b=0, p=-2\alpha.
$$
The solution: $u_{10}=3\alpha (01)$, $u_7=3\alpha (00)$. Similarly,
$u_{16}=3\beta (10)$, $u_{13}=3\beta (00)$. And the
other $u_j$ vanish. Having substituted them into (31)-(32) we get
$$
\left.\begin{matrix}
(31) &-3\alpha (00)-3\beta (00)=0\cr
(32) &6\alpha (10)-3\beta (10)=a\cr
\end{matrix}\right \}\Longrightarrow\alpha =\beta =0.
$$

4) $\lambda =0$, $\mu\geq 2$. From (66) we get $u_4=a(01)$,
$u_{10}=u_{16}=u_{20}=0$. This gives
$$
\left.\begin{matrix}
(63) &-\mu a=0\cr
(64) &-a+p=0\cr
 \end{matrix}\right \}\Longrightarrow u_4=u'=0.
$$
Moreover, equation (31) yields $u''=0$.

5) $\lambda =0$, $\mu =1$. From (66) we get $u_4=a(01)$, $u_{10}=3b(01)$,
$u_{16}=3c(01)$, $u_{20}=d(01)$.
Having substituted $u_4$, $ u_{10}$, $ u_{16}$, $ u_2$ as well as
$$
\begin{matrix}
 u_3=&X_+u_4=a(00), \cr
 u_9=&\frac{2}{3}X_+u_{10}+u'=(2b+p)(00), \cr
 u_7=&(b-p)(00), \cr
 u_{19}=&d(00), \cr
 u_{15}=&(2c+q)(00), \cr
 u_{13}=&(c-q)(00)
\end{matrix}
$$
directly into
(28)--(37) we get
$$
\begin{matrix}
(28) &u_3=0\Longrightarrow a=0\cr
(35) & u_{19}=0\Longrightarrow d=0\cr
(29) & x''_-u_9\Longrightarrow p=-2b, u_9=0, u_7=3b(00)\cr
(31) & u_7+u_{13}=0\Longrightarrow u_{13}=-3b(00), c-q=-3b\cr
(32) & x''_-u_{15}-u_{16}=0\Longrightarrow 2c+q-3c=0\Longrightarrow q-c=0
\Longrightarrow b=0.
\end{matrix}
$$
The solution is $u_{16}=3c(01)$, $u_{15}=3c(01)$. Having substituted it
into (42)--(44)
we get
$$
\left.\begin{matrix}
(42) &m\cdot 3c(00)=0\cr
(44) &(2m+1)\cdot 3c(01)=0\cr
\end{matrix}\right \}
\Longrightarrow f=0.
$$
There are no nonzero solutions.

6) $\lambda =\mu =0$. From (66) we get $u_4=u_{10}=u_{16}=u_{20}=0$,
$u'=u_9=-u_7=p(00)$, $u''=u_{15}=-u_{13}=q(00)$. From (31) we deduce that
$u_7+u_{13}=0\Longrightarrow q=-p$. Having substituted this
into (39), (42), (46) we get
$$
\begin{matrix}
(39) & lp=0\cr
(42) &(2l-1)p-mp(00)=0\cr
(46)&-mp=0.
\end{matrix}
$$
No nonzero solutions.

There remains the case $s=0$. In this case $u_4=a(00)$,
$u_{10}=b(00)$, $u_{16}=c(00)$, $u_{20}=d(00)$; so equation (29),
(32), (33) and (35) imply $u_4=0$, $u_{16}=0$, $u_{10}=0$, $u_{20}=0$,
respectively. No nonzero solutions. We have considered all the
cases.

\section{The general case $(n>2)$}

Recall that $V$ and $W$ are $\fgl(n)$-modules;
$$
\begin{array}{l}
 I(V^*)=\Kee[\partial_1,
\dots , \partial_n]\otimes V^* \\
I(V^*, W^*)=\Kee[\partial'_1, \dots ,
\partial'_n]\otimes V^*\otimes \Kee[\partial''_1, \dots , \partial''_n]
\otimes W^*\end{array}
$$
(the tensor product of the vector spaces, but not modules). Let $e_1,
e_2, \dots , e_n$ be a~basis in $\Kee^n$. Let $E=\Span(e_{i_1},
e_{i_2}, \dots , e_{i_j})\subset \Kee^n$. Denote
by $\fgl(E)\subset \fgl(n)$ the Lie algebra of the operators that
preserves $e_i\not\in E$.

Let $\cL_E=\Kee[\partial_{i_1}, \dots ,
\partial_{i_j}]$ be the subalgebra of $\cL$. Set
$$
I_E(V^*)=\Kee[\partial_{i_1}, \dots , \partial_{i_j}]\otimes V^*
,\quad I_E(V^*, W^*)=I_E(V^*)\otimes I_E(W^*).
$$

As $\fgl(E)$-modules, $V^*$ and $W^*$ split into the direct sum of irreducible
submodules:
$$
V^*=\oplus_\alpha V^*_{\alpha},\quad W^*=\oplus_{\beta} W^*_{\beta}.
$$
Hence, $V^*\otimes W^*=\oplus _{\alpha , \beta}V_{\alpha}\otimes V_{\beta}$.
This implies the decomposition
$$
I_E(V^*, W^*)=\oplus_{\alpha, \beta}
I(V^*_{\alpha}\otimes I(W^*_{\beta}).
$$

Denote by $\Pi_E$ or $\Pi_{i, i_2\dots , i_j}$ the natural projection
of $I(V^*, W^*)$ onto $I_E(V^*, W^*)$: we replace $\partial_i'$ and
$\partial_i''$ for $i\not\in (i_1, \dots , i_j)$ with zeros.

This projection commutes with the $\cL_E$-action, and therefore it
sends $\fgl(n)$-highest (hence,
$\fgl(E)$-highest) singular (with respect to $\cL$; hence, with respect to
$\cL_E$) vectors into the highest singular vectors.

By a~natural basis in $I(V^*, W^*)$ we will mean a~basis consisting of
the elements
$$
P(\partial'_1, \dots , \partial''_n)v\otimes w,
$$
where $P$ is a~monomial and $v$, $ w$ are elements of the Gelfand--Tsetlin
bases of $V^*$ and $W^*$, respectively (we will also denote such
elements by $P_1(\partial)v\otimes P_2(\partial)w$).

We will say that $f$ contains a~component $P(\partial'_1, \dots ,
\partial''_n)v\otimes w$ if the corresponding coordinate does not
vanish in the natural basis. The {\it type of the component} in this case
is the monomial $P(\partial'_1, \dots , \partial''_n)$.

Sometimes we will represent $f$ in the form
$$
f=\sum P_i(\partial'_1, \dots , \partial''_n)u_i,
$$
where $u_i\in V^*\otimes W^*$
and $P_i(\partial'_1, \dots , \partial''_n)u_i$ is the sum of all the
components of type $P_i$.

Recall that the weights of $\partial'_i$ and $\partial''_i$ are equal
to $(0, \dots , 0, -1, 0, \dots , 0)$ with $-1$ on the $i$-th place.
Among all the monomials that enter the decomposition of $f$, select the
monomial $P_0(\partial'_1, \dots , \partial''_n)$ with the least
(lexicographically) weight.

Let $f=P_0u_0+\sum P_iu_i$. Since $f$ is a~highest weight
vector, then
$$
(x_{\alpha}\partial_{\beta})f=0\text{ for $\alpha <\beta$ }
$$
but
$$
(x_{\alpha}\partial_{\beta})f=P_0\cdot
(x_{\alpha}\partial_{\beta})u_0\quad
\text{plus components of another type.}
$$

Hence, $(x_{\alpha}\partial_{\beta})u_0=0$, i.e., $u_0$ is a~highest
weight vector. Then $u_0=v_0\otimes w'+\dots +v'\otimes w_0$,
where $v_0, w_0$ are the highest weight vectors of $V^*$ and $W^*$,
respectively. Thus we have proved the following lemma.

\SSbegin{5.1}{Lemma}[On the highest component]\label{L5.1} Let $f\in I(V_1^*,
V_2^*)$ be a~highest weight vector (not necessarily singular). Let
$P_0(\partial_1, \dots , \partial_n)$ be one of the monomials with
lexicographically lowest weight among the monomials that enter the
decomposition of $f$. Then, $f$ contains the components $P_0v_0\otimes
w'$ and $P_0v'\otimes w_0$, where $v_0$ and $w_0$ are the highest
weight vectors of $V^*$ and $W^*$, respectively.
\end{Lemma}

The component $P_0v_0\otimes w'$ will be called
$V$-highest (or just the highest) while $P_0v'\otimes w_0$ will be called
$W$-highest one.

\section{Second order operators}

First, recall the list of the highest singular vectors or $n=1$, 2 found earlier.

{\bf $n=1$}.

a) $\lambda =\mu =(0)$, $\nu =(-2)$.
$$
f=\partial v\otimes w_0\partial w.
$$

b) $\lambda =(0)$, $\mu =(1)$, $\nu =(-1)$.
$$
f=\partial^2v\otimes w +\partial v\otimes \partial w.
$$

b${}'$) $\lambda =(1)$, $\mu =(0)$, $\nu =(-1)$.
$$
f=\partial v\otimes\partial w+v\otimes\partial^2 w.
$$

c) $\lambda =(0)$, $\mu =(m)$, $\nu =(m-2)$; ($m\neq 0, 1$).
$$
f=m\partial^2 v\otimes w+\partial v\otimes \partial w.
$$

c${}'$) $\lambda =(l)$, $\mu =(0)$, $\nu =(l-1)$.
$$
f=\partial v\otimes \partial w+lv\otimes \partial^2 w.
$$

d) $\lambda =(l)$, $\mu =(-l+1)$, $\nu =(-1)$; ($l\neq 0, -1$)
$$
f=(l-1)\partial^2v\otimes w+(2l-1)\partial v\otimes\partial w+
lv\otimes \partial^2 w.
$$

{\bf $n=2$}.
\begin{center}
\resizebox{\textwidth}{!}{
\begin{tabular}{|c|c|c|c|}
\hline
$n$&highest weights $\lambda $, $\mu$ and $\nu$&type of $\Pi_{1}f$&
type of $\Pi_{2}f$\cr
&of $V^{*}$, $W^{*}$ and $f$&&\cr
\hline
$1)$&$\lambda =\mu =(0, 0), \quad\nu =(-1, -1)$&-&-\cr
&$f=\partial_1v_0\otimes \partial_2w_0-\partial_2v_0\otimes \partial_1
w_0$&&\cr
\hline
$2)$&$\lambda =\mu =(0, 0), \quad\nu =(0, -2)$&-&a)\cr
&$f=\partial_2v_0\otimes\partial_2w_0$&&\cr
\hline
$3)$&$\lambda =(0, 0), \quad\mu =(1, -1), \quad\nu =(-1, -1)$&b)&b)\cr
&$2\partial_1^2v_0\otimes w+2\partial_1v_0\otimes\partial_1w_0+
 2\partial_1\partial_2 v_0\otimes w+\partial_1
 v_0\otimes\partial_2w_1+$&&\cr
&$\partial_2v_0\otimes \partial_1w_1+\partial_2^2v_0\otimes w_2+
 \partial_2v_0\otimes\partial_2w_2$&&\cr
\hline
$4)$&$\lambda =(0, -1), \quad\mu =(m, m), \quad\nu =(m-1, m-2)$&&c)\cr
&$f=m\partial_1\partial_2 v_0\otimes
w_0+\partial_1v_0\otimes\partial_2w_0+$&-&b) \text{ for } m=1\cr
&$m\partial_2^2v_1\otimes w_0+\partial_2v_1\otimes\partial_2w_0$&&a) \text{ for } m=0\cr
\hline
$4')$&$\lambda' =(l, l), \quad\mu =(0, -1), \quad\nu =(l-1, l-2)$&&c)\cr
&$f=\partial_2v_0\otimes\partial_1w_0+lv_0\otimes\partial_1\partial_2w_0+$&-&b') \text{ for } l=1\cr
&$\partial_2v_0\otimes\partial_2w_1+lv_0\otimes\partial_2^2 w_1$&&a) \text{ for } l=0\cr
\hline
$5)$&$\lambda =(0, -1), \quad\mu =(m+1, m), \quad\nu =(m-1,
m-1)$&c)&c)\cr
&$f=(m+1)\partial_1^2v_0\otimes w_0+\partial_1v_0\otimes\partial_1w_0+
 (m+1)\partial_1\partial_2v_0\otimes w_1+$&b) \text{ for } m=0&b)\cr
 &$(m+1)\partial_1\partial_2v_1
 \otimes w_0+\partial_1v_0\otimes\partial_2w_1+$&&\cr
&$\partial_2v_1\otimes\partial_1w_0+(m+1)\partial_2^2v_1\otimes w_1+
 \partial_2v_1\otimes\partial_2w_1$&a) \text{ for } m=-1&a)\cr
\hline
$5')$& is similar to 5)&&\cr
\hline
$6)$&$\lambda =(l, l), \quad\mu =(-l+1, -l), \quad\nu =(0, -1)$&&d)\cr
&$f=(l-1)\partial_1\partial_2v_0\otimes w_0+(l-1)\partial_1v_0\otimes
 \partial_2w_0+l\partial_2v_0\otimes \partial_1w_0+$&-&b) \text{ for }
 l=0\cr
&$lv_0\otimes\partial_1\partial_2w_0+(l-1)\partial_2^2v_0\otimes w_1+
(2l+1)\partial_2v_0\otimes \partial_2 w_1+lv_0\otimes \partial_2^2w_1$&&b') \text{ for } l=1\cr
\hline
$6')$& is similar to 6)&&\cr
\hline
\end{tabular}
}
\end{center}

A) First, consider the case when $\Pi_if=0$ for any $i$. For
$n=2$ this happens in case 1). The singular vector is of the form
$$
f=\partial_1v_0\otimes\partial_2w_0
$$ (here both $v_0$ and $w_0$ are
highest weight vectors of weight $\lambda =\mu =(0, 0)$ and the weight
of $f$ is equal to $\nu=(-1, -1)$).

In the general case the highest
component is of the form $\partial{i_0}v_0\otimes\partial_{j_0}w$
($i_0<j_0$). Hence, $\Pi_{i_0j_0}f\neq 0$ and is equal to the sum of
several vectors of the form 1), i.e.,
$$
\Pi_{i_0j_0}f=\sum_{\alpha,
\beta}a_{\alpha\beta}(\partial{i_0}v_{\alpha}
\otimes\partial_{j_0}w_{\beta}-\partial{j_0}v_{\alpha}\otimes\partial_{i_0}
w_{\beta}),
$$
where $v_{\alpha}\in V^*$, $w_{\beta}\in W^*$ are of weight $(\dots ,
0_{i_0},\dots , 0_{j_0}, \dots)$ with zeros on the $i$-th and $j$-th places.

Let $n=3$ for the moment. The following three cases are possible:

1) $i_0=1$, $j_0=2$. In this case the weight of $v_{\alpha}$ is equal
to $(0, 0, l)$ and the weight of $w_{\beta}$ is equal to $(0, 0, m)$,
where $l$ and $m$ do not depend on $\alpha$ and $\beta$ because the
sum of the coordinates of the weights are equal for all weight
vectors.

In particular, the weight of $v_0$ is equal to $(0, 0, l)$
and the weight of $w_0$ is equal to $(0, 0, m)$. The weight of $f$ is
equal to $(-1, -1, l+m)$.

Since $v_0$, $w_0$ and $f$ are highest weight vectors, it follows that
$$
l\leq 0, \qquad m\leq 0, \qquad l+m\leq -1.
$$
Since the multiplicity of the highest weight is equal to 1, it follows that
$\alpha$ and $\beta$ assume only one value, i.e.,
$$
\Pi_{12}f=a(\partial_1v_0\otimes\partial_2w_0-
\partial_2v_0\otimes\partial_1w_0.
$$
Observe that $\Pi_{12}f$ is not highest with
respect to $\fgl(3)$ because
$$
(x_1\partial_3)(\Pi_{12}f)=-a(\partial_3v_0\otimes\partial_2w_0-
\partial_2v_0\otimes\partial_3w_0)\neq 0,
$$
hence, $\Pi_{12}f\neq f$, i.e., either $\Pi_{13}f\neq 0$ or $\Pi_{23}f\neq 0$.

But $\Pi_{13}f$ and $\Pi_{23}f$ can be only of
the form 1) (or the sum of several vectors of the form 1)), and therefore 
$\nu_3=-1$, i.e., $l+m=-1$.

But $l\leq 0$, $m\leq 0$, hence, two cases are possible: $l=0$ or $m=-1$. In
both cases the operators exist:
$$
S_1(\varphi, w)=d\varphi\wedge dw, \quad S_1(w, \varphi)=dw\wedge d\varphi.
$$

Proof of uniqueness of the highest singular vector in each case.
Let
$$
f=a(\partial_1v_0\otimes\partial_2w_2-\partial_2v_0\otimes\partial_1w_0)+
\dots
$$
and
$$
g=b(\partial_1v_0\otimes\partial_2w_0-\partial_2v_0\otimes\partial_1w_0)+
\dots
$$
be highest singular vectors. Then, $bf-ag$ is  a~highest singular
vector but $\Pi_{12}(bf-ag)=0$.

Further in headings 2) and 3) we will
see that this case corresponds to other weights of $V^*$ and $W^*$.
Hence, $bf-ag=0$, i.e., $f$ and $g$ are proportional. Almost in all
the cases the uniqueness is also proved by this method.

Therefore, in what follows we will replace the proof with words ``the
uniqueness is proved routinely''. To apply the routine method, it
suffices to demonstrate that the highest component $P_0v_0\otimes w$
(or the $w$-highest component) is uniquely determined. In particular,
the weight of $w$ should be of multiplicity 1.

2) $i_0=1$, $j_0=3$. In this case the weight of $v_{\alpha}$ is equal
to $(0, l, 0)$ and the weight of $w_{\beta}$ is equal to $(0, m, 0)$;
the weight of $f$ is equal to $(-1, l+m, -1)$. Since the weight of
$f$ is a~highest one, then $l+m=-1$ and this means that either the
weight of $v_{\alpha}$ or the weight of $w_{\beta}$ is not highest
contradicting to Lemma \ref{L5.1}. Hence, there are no
singular vectors.

3) $i_0=2$, $j_0=3$. The weight of $v_{\alpha}$ is equal to $(l, 0,
0)$, the weight of $w_{\beta}$ is equal to $(m, 0, 0)$ and the weight
of $f$ is equal to $(l+m, 0, 0)$. Lemma \ref{L5.1} implies
that the weights $(l, 0, 0)$ and $(m, 0, 0)$ are highest ones, hence,
of multiplicity 1, and therefore 
$$
\Pi_{23}f=a(\partial_2v_0\otimes\partial_3w_0-
\partial_3v_0\otimes\partial_2w_0).
$$
Since $\Pi_{13}f=\Pi_{12}f=0$, it
follows that
$$
f=\Pi_{23}f=a(\partial_2v_0\otimes\partial_3w_0-
\partial_3v_0\otimes\partial_2w_0).
$$
The condition $(x^2_3\partial_1)f=0$ yields
$$
0=(x^2_3\partial_1)f=-2a(\partial_2v_0\otimes\partial_1((x_3\partial_1)w_0)-
\partial_1((x_3\partial_1)v_0)\otimes\partial_3w_0)
$$
implying $(x_3\partial_1)v_0=(x_3\partial_1)w_0=0$. This is true only 
for $l=m=0$. The corresponding operator exists. It is
$$
S_1(\varphi, \psi)=d\varphi\wedge d\psi.
$$
The case $n=3$ is considered completely.

Let us pass to the general case. First, let us prove that
$j_0=i_0+1$. Indeed, otherwise $\Pi_{i_0, i_0+1, j_0}f$ should be of
type 2) but there are no such vectors.

Let us show that the weight of $f$ is equal to $\nu=(0, \dots , 0, -1,
-1, \dots , -1)$ (with $(i_0-1)$-many 0's); the weight of $V^*$ is
$\lambda=(0, \dots 0, \underbrace{-1, \dots , -1}_{l})$; the weight of
$W^*$ is $\mu =(0, \dots , 0, \underbrace{-1, \dots , -1}_{m})$.

Indeed, for $i<i_0$ the vector $\Pi_{ii_0j_0}f$ is the sum of several
components of the form (3), hence, $\nu_i=0$.

Moreover, $\Pi_{ii_0j_0}f$ contains the highest component
$l_{i_0}v_0\otimes l_{j_0}w'$ implying $\lambda_i=0$. For $i>j_0$ the
vector $\Pi_{i_0j_0i}f$ is of the form 1), hence, $\nu_i=-1$. But
$\Pi_{i_0j_0i}f$ again contains the highest component implying either
$\lambda_i=0$ or $\lambda_i=-1$.

Similarly, the weight of $W^*$ is equal to $(0, \dots , 0, -1, \dots ,
-1)$. From the balance of the sum of weight coordinates it follows
that $l+m+2=n-i_0+1$. Moreover, $l\leq n-1$, $m\leq n-1$
$\lambda_{i_0}=\mu_{i_0}=0$ occurs.

In all these cases the invariant operators exist. These operators are
of the form
$$
S_1(w_1, w_2)=dw_1\wedge dw_2.
$$
The uniqueness is proved routinely, since the highest component
$l_{i_0}v_0\otimes l_{j_0}w$ is uniquely determined (namely,
$i_0=n-l-m-1$, $j_0=n-l-m$, and the weight of $w$ is equal to $(0,
\dots , 0, \underbrace{-1, \dots , -1}_{m}, \underbrace{0, \dots ,
0}_{l})$ and $W^*$ has no multiple weights). Thus, case A) is
considered completely.

B) Let now there exist an $i$ such that $\Pi_if\neq 0$. Let $i_1$ be
the least such index. Let $\Pi_{i_1}f=\sum_{\alpha}f_{\alpha}$ be the
sum of several vectors of types a)--d).

1) Let at least one of the summands be of type a). Then, 
$\nu_{i_1}=-2$, and therefore all of them are of type a). If
$i>i_1$, then $\Pi_{i_1i}f$ should be of type 5), that is $(m=-1)$,
implying $\nu_i=-2$. For $i<i_1$ we see that $\Pi_{ii_1}f$ is of type
4), i.e. $(m=0)$, or of type 2) implying either $\nu_i=0$ or
$\nu_i=1$.

Thus, the weight of $f$ is equal to $\nu=(0, \dots , 0,
\underbrace{-1, \dots , -1}_{p}, \underbrace{-2, \dots , -2}_{q})$,
where $q\geq 1$.

Now, let us show that if $\nu =(0, \dots , 0, -1, \dots , -1, -2,
\dots , -2)$, then $\lambda $ and $\mu$ are of the form $(0, \dots ,
0, -1, \dots , -1)$ with $l$ (resp. $m$) $-1$'s.

First, let $n=3$. The following cases are possible.

1) $\nu =(0, 0, 2)$. Observe that $\Pi_{12}f=0$ because in the list of
singular vectors for $n=2$ there is no vector $g$ such that
$\nu_1=\nu_2=0$ and $\Pi_1g=\Pi_2g=0$. Moreover, $\Pi_{13}f=\Pi_3f$
and $\Pi_{23}f=\Pi_3f$ because they are of type 2).

Therefore, the $V$-highest and $W$-highest components are of the form
$\partial_3v\otimes\partial_3w$, where the
weights of $v$ and $w$ are equal to
$(0, 0, 0)$ because $\Pi_{13}f$ and $\Pi_{23}f$ are of type 2).

2) $\nu =(0, -1, -2)$. Again $\Pi_{12}f=0$, $\Pi_{13}f$ is of type
2),
$\Pi_{23}f$ of type either 4) or 4${}'$) ($m=0$). The highest
component is of the form
$$
\partial_2v_0\otimes\partial_3w\text{ or }\partial_3v_0\otimes\partial_2w.
$$

In the first case $\Pi_{13}f=\partial_3v\otimes\partial_3w$ is of the
form 2) and
$$
\Pi_{23}f=\partial_2v_0\otimes\partial_3w_0+\partial_3v\otimes\partial_3w_0
$$
is of type 4, hence, the weight of $v$ is equal to $(0, -1, 0)$ and
the weight of $v_0$ is equal to $(0, 0, -1)$; in the second case
$$
\Pi_{23}f=\partial_3v_0\otimes\partial_2w_0+
\partial_3v_0\otimes\partial_3w
$$ is of type 4${}'$) and the weight of
$v_0$ is equal to $(0, 0, 0)$.

Similarly, the weight of $w_0$ in the first case is equal to $(0, 0, 0)$ and in
the second one
is equal to $(0, 0, -1)$.

3) $\nu=(-1, -1, -2)$. The vectors $\Pi_{13}f$ and $\Pi_{23}f$ are of
type 4) or 4${}'$). There are three possibilities:

3.1) $\Pi_{13}f$ and $\Pi_{23}f$ are of type 4). Then, $\Pi_3f$ consists
of the components of the form $\partial_3v\otimes\partial_3w$
where the weight of $v$ is equal to $(-1, -1, 0)$ and the weight of $w$ is
equal to $(0, 0, 0)$. The
image under $\Pi_{13}$ consists of components
$$
\partial_1v_0\otimes\partial_3w+\partial_3v\otimes\partial_3w,
$$
where $v_0=(x_1\partial_3)v$. Therefore, the weight of $v_0$ is
equal to $(0, -1, -1)$.

If we would have proven that
$\partial_1v_0\otimes\partial_3w$ is the highest component this would
have implied that $\lambda =(0, -1, -1)$, $\mu =(0, 0, 0)$.

But the component
$\Pi_{12}f=\partial_1v\otimes\partial_2w+\dots , $ of type 1) is also possible
and then the weight of $v$ could be equal to $(0, 0, -2)$ while the weight of
$w$ to $(0, 0, 0)$.

But observe
that
$$
(x_1\partial _2)(\Pi_3f)=a\partial_3(x_1\partial_2)v\otimes\partial_3w
$$
cannot cancel with other components of $(x_1\partial_2)f$. Hence,
$(x_1\partial_2)(\Pi_3f)=0$, i.e., $(x_1\partial_2)v=0$ which is
impossible if the weight of $V^*$ is equal to $(0, 0, -2)$. Therefore,
$\Pi_{12}f=0$ and the highest weight
of $V^*$ is equal to $(0, -1, -1)$.

3.2) $\Pi_{13}f$ and $\Pi_{23}f$ are of type 4${}'$ is treated similarly.

3.3) Both types 4) and $4')$ are encountered. Then, $\Pi_3f$ consists of
components of the form $\partial_3v\otimes\partial_3w$, where either
the weight of $v$ is equal to $(0, -1, 0)$ and the weight of $w$ is equal to
$(-1, 0, 0)$ or the other way round.

The components of type
$\partial_1\partial_3$ are vectors $\partial_1v\otimes\partial_3w$,
where the weight of $v$ is equal to $(0, -1, 0)$ and the weight of $w$
is equal to $(0, 0, -1)$ (if the corresponding summand $\Pi_{13}f$ is
of type 4)) and $\partial_3v\otimes\partial_1w$, where the weight of
$v$ is equal to $(0, 0, -1)$ and the weight of $w$ is equal to $(0,
-1, 0)$ (if the corresponding summand is of type $4')$).

The components
of $\Pi_{12}f$ (if any) are of the form
$$
\text{ $\partial_1v\otimes\partial_2w$ and
$\partial_2v\otimes\partial_1w$},
$$
where the weights of $v$ and $w$ are equal to $(0, 0, -1)$ because
$\Pi_{12}f$ is of type 1). Among the listed components there are
highest ones, hence, the highest weights of $V^*$ and $W^*$ are equal
to $(0, 0, -1)$.

4) $\nu =(0, -2, -2)$. The $V^*$-highest and
$W^*$-highest components are 
$\partial_2v\otimes\partial_2w$. The vector $\Pi_{12}f$ is of type 2)
$\Pi_{23}f$ of type 5) ($m=-1$), hence, the weights of $v$ and $w$
are equal to $(0, 0, -1)$.

5) $\nu =(-1, -2, -2)$. The vector
$\Pi_{23}f$ is of type 5) (if $m=-1$), $\Pi_{12}f$ of type 4) or
$4')$ (if $m=0$).

In the first case the weight of $v$ is equal to $(-1, 0,
-1)$, the weight of $w$ is equal to $(0, 0, -1)$ while in the second
case it is the other way round.

In the first case
$$
\Pi_{12}f=\partial_2v_0\otimes\partial_2w+\partial_1v\otimes\partial_2w,
$$
where $v_0=(x_1\partial_2)v$, hence, the weight
of $v_0$ is equal to $(0, -1, -1)$ and the weight of $w$ is equal to
$(0, 0, -1)$.

In the second case the
weights are transposed.

6) The $\nu =(-2, -2, -2)$. The $V^*$ -highest and $W^*$-highest components
are of the form $\partial_1v\otimes\partial_1w$; the vectors
$\Pi_{12}f$ and $\Pi_{13}f$ are of type 5), hence, the weights of $v$
and $w$ are equal to $(0, -1, -1)$.

Let now $n\geq 4$. Let
$P(\partial_{i_0}\partial_{j_0})v\otimes w$ be the highest component,
the weight of $v$ be equal to $(\lambda_1, \dots , \lambda_n)$. Then, 
$\Pi_{i_0j_0i}f$ is of one of the types indicated, hence, either
$\lambda_i=0$ or $\lambda_i=-1$.

Similarly, the highest weight of $W^*$
is equal to $(0, 0, \dots , -1)$. Thus,
$$
\lambda =(0\dots , 0, \underbrace{-1, \dots , -1}_{l}, \;
\mu =(0, \dots , 0\underbrace{-1, \dots , -1}_{m}, \;
\nu =(0, \dots , \underbrace{-1, \dots , -1}_{p}, \underbrace{-2, \dots , -2}_{q}.
$$
The corresponding operators should be conjugate to those we will
consider in the next heading. Therefore, we will not prove neither
their existence
nor their uniqueness.

2. The second case: $\Pi_{i_1}f$ contains a~summand of the form
$$
\partial^2_{i_1}v\otimes w+\partial_{i_1}v\otimes\partial_{i_1}w
$$
but does not contain any summands of type 2). Then, $\nu =-1$ and all
the summands are of the form indicated; the $i_1$-th coordinate of the
weight $v$ vanishes while that of $w$ is equal to 1. $\Pi_{i_1i}f$
for $i>i_1$ is either of type 3) or of type 5) $(m=0)$ , hence,
$\nu_i=1$ and $\Pi_if$ is also of type b). $\Pi_if=0$ for $i<i_1$
and,
therefore, $\Pi_{ii_1}f$ is either of type 4) ($m=1$) or 6) ($l=0$)
implying $\nu_i=0$.

Thus, $\nu =(0, \dots , 0, \underbrace{-1, \dots , -1}_{r})$.
If $i, j<i_1$, then $\Pi_{ij}f=0$ (the list of singular vectors
for $n=2$ does not contain any vector $g$
such that $\nu_1=\nu_2=0$ and $\Pi_1g=\Pi_2g=0$).

Therefore, the $V$-highest
and $W$-highest components are $\partial_1\partial_{i_1}v_0\otimes w$ and
$\partial_1\partial_{i_1}v\otimes w_0$, respectively (perhaps, $i_1=1$).

Denote the weights of $v_0$, $v$, $w_0$, $w$ by
$$
(\lambda_1, \dots , \lambda_n), \quad
(\lambda'_1, \dots , \lambda'_n), \quad (\mu_1, \dots , \mu_n), \quad
(\mu'_1, \dots , \mu'_n).
$$
Since $\Pi_{ii_1}f$ is either of the form 4) or 6), it
follows that $\lambda_1=\lambda'_1=0$ and $\mu_1=\mu'_1=1$ (if $i_1=1$, then
this is also true because $\Pi_{i_1}f$ is of type b)).

But
$\lambda_1=\max\{\lambda_i, \lambda'_i\}$ and $\mu_1=\max\{\mu_i, \mu'_i\}$,
hence, $\lambda_i\leq 0$, $\mu_i\leq 1$, $\lambda'_i\leq 0$, $\mu'_1\leq 1$.
Observe that
$$
\lambda'_i+\mu_i=\left\{\begin{matrix}
\nu_i, \text{ for }i\neq1, i\neq i_1\cr
\nu_i+1, \text{ for }i=1\text{ or }ii=i_1, \; i_1\neq 1\cr
\nu_i+2, &\text{ for }ii=i_1=1\cr
\end{matrix}\right .
$$
But either $\nu_i=0$ or $\nu_i=1$, hence, $\lambda'_i+\mu_i\geq -1$.

But $\lambda'_i\leq 0$, hence, $-1\leq\mu_i\leq 1$. Thus,
$(\underbrace{-1, \dots , -1}_{p}, 0, \dots , 0, \underbrace{-1,
\dots , -1}_{q})$ is the highest weight of $W^*$, and therefore
the highest weight of $(W^+)^*$ is equal to
$$
(0, \dots , -1, \dots , -1, -2, \dots , -2). $$
So every operator
encountered in this heading is conjugate to the operator from the
previous heading.

In particular, this implies that the highest weight of $V^*$ is $(0,
\dots , 0, \underbrace{-1, \dots , -1}_{l})$.

\SSbegin{5.2}{Statement}\label{S5.2} {\em 1)} $l\leq n-1$; {\em 2)} $ p\geq 1$;
{\em 3)} $r\geq 1$; {\em 4)} $ l+1\geq p$; {\em 5)} $l-p+q+2=r$.
\end{Statement}

\begin{proof} 1) It follows from the fact that $\lambda_1=0$.

2) It follows from the fact that $\mu_1=1$.

3) It follows from the fact that $\nu_{i_1}=-1$.

4) Let us consider the $W$-highest component
$\partial_1\partial_{i_1}v\otimes w_0$. We have $\mu_i=1$ and
$\lambda'_i+\mu_i\leq 0$ for $i=2, 3, \dots , p$; hence, $\lambda'_i=-1$.
Therefore, the number of $-1$'s among the coordinates of the weight of
$v$ is equal to
$p-1$, i.e., $l\geq p-1$.

5) It is verified by the direct comparison of the sums of coordinates of the
weights.
\end{proof}

\SSbegin{5.3}{Statement}\label{S5.3} For any collection of $l$, $p$, $q$, $r$ satisfying the
conditions of Statement \ref{S5.2}, there exists an invariant operator and this
operator is unique.
\end{Statement}

\begin{proof} {\bf Existence}. Observe that in the tensor product of
$\Lambda^{l+1}(\Kee^n)$ and the representation with highest weight
$(\underbrace{1, \dots , 1}_{p}, 0, \dots , 0, \underbrace{-1, \dots ,
-1}_{q})$ contains an irreducible component isomorphic to
$\Lambda^{r-1}(k^n)$. Therefore, the operator $S_1^{*-1}(w, s)=
dZ(dw, s)$ exists.

{\bf Uniqueness}. Observe that the
weights of $V^*$ are multiplicity free hence, the $W$-highest component
$\partial_1\partial_{i_1}v\otimes w_0$ is uniquely
determined. Apply the standard method.
\end{proof}

3) $\Pi_{i_1}f$ contains the summands of the form
$$
m\partial^2_{i_1}v\otimes w +\partial_{i_1}v\otimes\partial_{i_1}w
(m\neq 0, 1).
$$
Then, $\nu_{i_1}=m-2$, and therefore all the summands
in $\pi_{i_1}f$ are of this form.

The vector $\Pi_{i_1i}f$ for $i>i_1$ is of
type 5) implying $\nu_i=m-2$.

The vector $\Pi_{ii_1}f$ for $i<i_1$ is of type
4), hence, $\nu_i=m-1$ and the component
$\partial_i\partial_{i_1}v\otimes w$ does not vanish.

We have $\Pi_{ij}f=0$ for $i, j<i_1$ because for $n=2$ there is no
such a~vector. Hence, the $V$-highest and $W$-highest components are
of the form $\partial_1\partial_{i_1}v\otimes w$. The vector
$\Pi_{ii_1}f$ is of type 4), and therefore contains components
$$
\partial_1\partial_{i_1}v\otimes w+\partial^2_{i_1}v'\otimes w,
$$
where
$v=(x_1\partial_{i_1})v'$. Since $\Pi_{ii_1}f$ for $i>i_1$ is of type 5) and
$i<i_1$ of type 4) for $v$, it follows that the weight of $v$ is
equal to $(-1, \dots , -1, 0, -1, \dots , -1)$ while the weight
of $w$ is equal to $(\underbrace{m, \dots , m}_{i_1}, m-1\dots , m-1)$.

Since $v=(x_1\partial_{i_1})v'$, it follows that the weight of $v$ is equal to
$(0, -1, \dots , -1)$.
Thus,
$$
\lambda=(0, -1, \dots , -1), \quad \mu =(\underbrace{m, \dots , m}_{i_1},
m-1, \dots , m-1),
$$
$$
\nu =(\underbrace{m-1, \dots , m-1}_{i_1-1}, m-2, \dots , m-2).
$$
The corresponding operator is $S_2(w, s)=P_4(dw, s)$ and its
uniqueness is proved routinely.

4) The vector $\Pi_{i_1}f$ contains a~component of type d), $\nu_{i_1}=-1$.
The vector $\Pi_{ii_1}f$ vanishes for $i>i_1$; hence, $i=n$.

The vector $\Pi_{in}f$ is of type 6) for $i<n$, hence, $\nu_i=0$.
Thus, $\nu =(0, \dots , 0, -1)$. This already shows that the
operator conjugate to this one is of different type. But since all the other
cases are already considered, it follows that in this case all the operators are
conjugate to the ones already considered.

We have considered all the
possibilities for the singular vectors of degree 2.

\section{Third order operators}

First, let us recall the list of singular highest weight vectors of degree 3 for
$n=2$.

1) $\lambda =\mu =(0, -1)$, $\nu =(-2, -3)$.
$$

\begin{array}{l}
f=\partial_1\partial_2v_0\otimes\partial_1 w_0 -\partial_1v_0\otimes
\partial_1\partial_2w_0+\partial^2_2v_1\otimes\partial_1 w_0+\\
\partial_1\partial_2v_0\otimes\partial_2 w_1-\partial_1v_0\otimes
\partial^2_2w_1 -\partial_2v_1\otimes\partial_1\partial_2 w_1-
\partial_2v_1\otimes\partial^2_2 w_1,
\end{array}
$$
where $v_0=X_+v_1$ and $w_0=X_+w_1$; the weights of $v_0$ and $w_0$
are equal to $(0, -1)$; the weights of $v_1$ and $w_1$
are equal to $(-1, 0)$.

2) $\lambda =\nu =(0, -1)$, $\mu =(2, 1)$. The singular vector is of
the form
$$

\begin{array}{l}
f=2\partial^2_1\partial_2v_0\otimes w_0 +\partial_1\partial_2v_0\otimes
\partial_1 w_0 +2\partial^2_1v_0\otimes\partial_2 w_0 +\\
\partial_1v_0\otimes\partial_1\partial_2 w_0+2\partial_1\partial^2_2v_0
\otimes w_1 +3\partial_1\partial_2v_0\otimes\partial_2 w_1+
\partial_1v_0\otimes\partial^2_2 w_1+\\
\partial_1\partial^2_2v_1\otimes w_0+\partial^2_2v_1\otimes\partial_1 w_0 +
2\partial_1\partial_2v_1\otimes\partial_2 w_0-
\partial_2v_1\otimes\partial_1\partial_2 w_0 +2\partial^3_2v_1\otimes
w_1+\\
3\partial^2_2v_1\otimes\partial_2 w_2 +
\partial_2v_1\otimes\partial^2_2 w_1,
\end{array}
$$
where $v_0=X_+v_1$ and $w_0=X_+w_1$; the weights of $v_0$, $v_1$,
$w_0$, $w_1$ are equal to $(0, -1)$, $(-1, 0)$, $(2, 1)$, $(1, 2)$,
respectively.

$\lambda =(2, 1)$, $\mu =\nu =(0, -1)$. This case is similar to 2).
Consider two subcases:

a) We have $\Pi_{ij}=0$ for any $i, j$. Then, there exists a~space
$E=\langle e_{i_1}, e_{i_2}, e_{i_3}\rangle$ such that
$$

\begin{array}{l}
\Pi_if=\\
\partial'_1\partial'_2\partial'_3u_1+\partial'_1\partial'_2\partial''_3u_2+
\partial'_1\partial''_2\partial'_3u_3+\partial''_1\partial'_2\partial'_3u_4+\\
\partial'_1\partial''_2\partial''_3u_5+\partial''_1\partial'_2\partial''_3u_6+
\partial''_1\partial''_2\partial'_3u_7+\partial''_1\partial''_2\partial''_3u_8
\end{array}
$$
But this vector cannot be a~highest weight one. Indeed,
$$

\begin{array}{l}
0=(x_1\partial_2)f=\\
-\partial'{}^2_2\partial'_3u_1-\partial'{}^2_2\partial_3u_2-
\partial'_2\partial''_2\partial'_3(u_3+u_4)-
\partial'_2\partial''_2\partial''_3(u_5+u_6)-\partial''{}^2_2\partial'_3u_7-
\partial''{}^2_2\partial''_3u_8
\end{array}$$
+ components of the same form as in
$f$.

Hence, $u_1=u_2=u_7=u_8=u_3+u_4=u_5+u_6$ and
$$

\begin{array}{l}
0=(x_2\partial_3)f=-\partial'_1\partial'_3\partial''_3u_3-\\
\partial''_1\partial'{}^2_3u_4-\partial'_1\partial''{}^2_3u_5-
\partial''_1\partial'_3\partial''_3u_6+ \\
\text{components of the same form as in $f$;}
\end{array}
$$
hence, $u_3=u_4=u_5=u_6=0$, i.e., $f=0$.

B) There exist $i, j$ such that $\Pi_{ij}f\neq 0$. Then, $\Pi_{ij}f$ is
the sum of several vectors of types 1), 2)
and 2'). Moreover, if at least one of the summands is of type 1), then
$\Pi_{i'j}f$ and, therefore, all the summands are of the same type.

Moreover, since $\nu_j=-3$ and $\Pi_jf\neq 0$, then $\Pi_{i{'}j}f$
is also of type 1) for $i'<j$. For $j'>j$ the vector $\Pi_{jj{'}}f$ vanishes
implying that $j=n$.

Since $\Pi_{1n}f$ is of type 1), the weights of the vectors $v$ and $w$ that
enter $\partial^2_nv\otimes\partial_nw$ are
equal to $(-1, \dots , -1, 0)$ and the weight of $(-2, \dots , -2, -3)$ is equal to
$f$.

So far, the highest of the components found is
$\partial_1\partial_nv'\otimes\partial_1w'$. It enters the vector
$$
\Pi_{1n}f=\partial_1\partial_nv'\otimes\partial_1w'+ \dots +
\partial^2_nv\otimes\partial_nw,
$$
where $v'=(x_1\partial_n)v$ and $w'=(x_1\partial_n)w$, and therefore 
the weights of $v'$ and $w'$ are equal to $(0, -1, \dots , -1)$.

But $\Pi_{1j}=0$ for $j<n$ because $\nu_1=\nu_j=-2$ and no such
operator exists for $n=2$. Therefore, the $V$-highest and $W$-highest
components are of the form $\partial_1\partial_nv'\otimes\partial_1w'$
implying that $\lambda =\mu =(0, -1, \dots , -1)$. Such an operator
exists:
$$
T(w_1, w_2)=P_4(dw_1, dw_2).
$$

If $\Pi_{ij}f$ does contain a~summand of type 2) or $2')$, then
$\nu_j=-1$ and all the other summands are of the same form.

We similarly prove that $j=n$ and the weight of $f$ is equal to $(0,
\dots, 0, -1)$. But either the highest weight of $V^*$ or the highest
weight of $W^*$ is of different form because $\Pi_{in}f$ is either of
type 2 or $2')$. Therefore, in this case the operators obtained are
conjugate to the operators of different types, i.e., to $T(w_1, w_2)$.
The proof is completed.

\section{Operators in the spaces of twisted forms}

Recall that the space twisted of $p$-forms with twist $l$ is $T(l+1,
\dots , l+1, l, \dots, l)$ with $p$-many $l+1$'s. Denote the space
$\Omega^p_{l}:=\Omega^p\otimes \Vol^l$. In particular, $p$-forms with
twist $0$ are the usual differential forms, $p$-forms with twist $-1$
are polyvector fields, any volume form on the $n$-dimensional manifold
$M$ can be considered as a~$0$-forms with twist $1$ or as an $n$-form
with twist 0.

Having selected a~nondegenerate volume form $\delta$ on $M$, any
element from $\Omega^p_{l}$ can be represented either in the form
$$
s=\omega\cdot \delta^{l}, \text{ where $\omega\in\Omega^p$, }
$$
thanks to the existence of the nonzero $0$-order operator $Z: \Omega^p\otimes
\Vol^l\tto \Omega^p_{l}$, or, if alternatively, in the form
$$
s=\xi\cdot \delta^{l+1}, \text{ where $\xi\in\Omega^p_{-1}•$ is a
polyvector of degree $n-p$}.
$$

Let $x=(x_1, \dots , x_n)$ be local coordinates in a~neighborhood of a
point $P\in M$. Select $\delta=\vvol$. Obviously, volume preserving
diffeomorphisms transform twisted forms by the same formulas as the
usual forms: they consider $\Vol^l$ as the space of functions.

\ssec{Zero order differential operators} Let the weight of the
$\fgl(n)$-module $V$ be
$$
\lambda= (l+1, \dots , l+1, l, \dots, l)$$ with $p$-many $l+1$'s;
that of the $\fgl(n)$-module $W$ be
$$\mu= (m+1, \dots , m+1, m,
\dots, m)$$ with $q$-many $m+1$'s. Then, $V\otimes W$ splits into
the direct sum of irreducible modules of which exactly one has the
highest weight of the same form:\footnotesize
$$
\nu=\begin{cases}(l+m+1, \dots , l+m+1, l+m, \dots, l+m)&\text{with
$p+q$-many $l+m+1$'s }\cr &\text{if $p+q\leq n$}\cr (l+m+2, \dots ,
l+m+2, l+m+1, \dots, l+m+1)&\text{with $p+q-n$-many $l+m+2$'s}\cr
&\text{if $p+q\geq n$}\end{cases}
$$\normalsize
The corresponding invariant operators are the exterior product of
twisted differential forms and the exterior product of twisted
polyvector fields; in both cases, the ``twists'' are considered as
coefficients:
$$
Z_{1}(\omega_1\cdot \delta^{l}, \omega_2\cdot
\delta^{m})=\omega_1\wedge\omega_2\cdot \delta^{l+m}
$$
$$
Z_{1}(\omega_1\cdot \delta^{l}, \omega_2\cdot
\delta^{m})=\omega_1\wedge\omega_2\cdot \delta^{l+m}\; \text{for $p+q\leq n$};
$$
$$
Z_{2}(\xi_1\cdot \delta^{l+1}, \xi_2\cdot
\delta^{m+1})=\xi_1\wedge\xi_2\cdot \delta^{l+m+2}\; \text{for $p+q\geq n$}.
$$
From multiplicity-free occurrence of the target space, it follows
that $Z_{1}$ and $Z_{2}$ are proportional if $p+q=n$.

I hope the reader can forgive me that here I skip verification of
invariance of these operators. Actually, to prove it correctly with
all details takes some space and arguments.

\ssec{First order differential operators} On the spaces of twisted
forms, there are only the following invariant bilinear operators:

1) $P_6: \Omega^p_{l}\times \Omega^q_{m}\tto \Omega^{p+q+1}_{l+m}$
for $p+q\leq n-2$ and $(l, m)\neq (0, 0)$;

2) $P_7: \Omega^p_{l}\times \Omega^q_{m}\tto \Omega^{p+q+1-n}_{l+m+1}$
for $p+q\geq n-1$ and $(l, m)\neq (0, 0), (0, -1), (-1, 0)$;

in the exceptional cases we have

3) $P_1(\omega_{1},
\omega_{2})=ad\omega_{1}\wedge\omega_{2}+b\omega_{1}\wedge
d\omega_{2}$ for $p+q\leq n-2$ and $(l, m)=(0, 0)$;

4) $P_1(\omega_{1}, \omega_{2})=aZ_{2}(d\omega_{1},
\omega_{2})+bZ_{2}(\omega_{1}, d\omega_{2})$ for $p+q\geq n-1$ and
$(l, m)=(0, 0)$;

5) {\small $P(\omega_{1}, \omega_{2}\delta^{-1})=aZ_{2}(d\omega_{1},
\omega_{2}\delta^{-1})+bdZ_{2}(\omega_{1}, \omega_{2}\delta^{-1})$ for $p+q\geq n-1$ and $(l, m)=(0, -1)$}; for $(l, m)=(-1, 0)$ {\it mutatis mutandis};

additionally, if the result is not a~twisted form, there exist the
following operators:

6) $P_1(\omega, s)=Z(d\omega, s)$ for $l=0$;

7) $P_1(s, \omega)=Z(s, d\omega)$ for $m=0$.

I hope that the reader forgives me not retyping one more page of
formulas in order to prove the more or less obvious.

Let us prove that the above list exhausts all the operators. To this
end, let us prove that the highest singular vectors exist only in
these cases.

{\bf For $n=1$} the singular vectors are:

\noindent 1) $\lambda =(l)$, $\mu =(m)$, $\nu
=(l+m-1)$ and $(l, m)\neq (0, 0)$:

$f=m\partial v\otimes w- lv\otimes\partial w$,

\noindent 2) if $(l, m)=(0, 0)$, then
$f=a\partial v\otimes w+bv\otimes\partial w$ for any $a, b$.

{\bf For $n=2$} the singular vectors are:

\noindent 1) $l=m=0$. In particular, $\lambda =(l, l-1)$, $\mu =(m, m-1)$, $\nu
=(l+m-1, l+m-2)$,
$$
\begin{array}{l}
f=(l+m-1)(m\partial _{1}v_0\otimes w_0-
lv_0\otimes\partial _{1}w_0)+ l(m\partial_2v_0\otimes w_1-
(l-1)v_0\otimes\partial_2w_1)+ \\
m((m-1)\partial_2v_1\otimes w_0-
lv_1\otimes\partial_2w_0),\end{array}
$$

$f=a(\partial_1v_0\otimes w_0 +
\partial_2v_1\otimes w_0)+ b(v_0\otimes\partial_2w_0 +
v_0\otimes\partial_2w_1)$,

$f=(a+b)v_0\otimes\partial_1w_0
-a\partial_2v_0\otimes w_1)+ bv_0\otimes\partial_2w_1 +
a\partial_2v_1\otimes w_0 + av_1\otimes\partial w_0$

\noindent 2) $l=m=0$. In particular, $\lambda =(l, l)$, $\mu =(m, m-1)$, $\nu
=(l+m-1, l+m-1)$,

$f=m\partial_1v_0\otimes
w_0-lv_0\otimes\partial_1w_0+ m\partial_2v_0\otimes
w_1-lv_0\otimes\partial_2w_1$,

$f=a\partial_1v_0\otimes w_0+
bv_0\otimes\partial_1w_0+ a\partial_2v_0\otimes w_1+
bv_0\otimes\partial_2w_1$

\noindent $2'$) $\lambda =(l, l-1)$, $\mu=(m, m)$, $\nu=(l+m-1), l+m-1)$ is similar.

\noindent 3) $l=m=0$. In particular, $\lambda =(l, l)$, $\mu =(m, m)$, $\nu =(l+m,
l+m-1)$,

$f=m\partial_2v_0\otimes w_0-lv_0\otimes\partial_2w_0$,

$f=a\partial_2v_0\otimes w_0+bv_0\otimes\partial_2w_0$

\noindent 4) $\lambda =(0, 0)$, $\mu =(m, m-1)$, $\nu =(m, m-2)$,

$f=\partial_2v_0\otimes w_0$.

\noindent $4')$ $\lambda =(l, l-1)$, $\mu =(0, 0)$, $\nu =(l, l-2)$,

$f=v_0\otimes\partial_2w_0$.

{\bf $n\in \Zee_{+}•$}. Let us prove that $\Pi_nf\neq 0$. Indeed,
$\Pi_{in}f$ is of one of the types 1)--4), hence, $\Pi_nf\neq 0$. The
vector $\Pi_{1n}f$ is also of one of these types, hence,
$\nu_1-\nu_n\leq 2$.

Consider the case $\nu_1-\nu_n=2$. Then, $\Pi_{1n}f$ is the sum of the
summands of type 4) or 4${}'$). Let at least one of the summands,
$\partial_nv\otimes w$, be for example of type 4).

Let the weights of $v$ and $w$ be equal to $(\lambda_1, \dots ,
\lambda_n)$ and $(\mu_1, \dots , \mu_n)$, respectively. Then, 
$\lambda_n=0$, $\mu_1=m$, $\mu_n=m-1$. Let $f$ be the weight of
\[
    (
        t, \dotsc, t,
        \underbrace{t-1, \dotsc, t-1}_{a},
        \underbrace{t-2, \dotsc, t-2}_{b}
    ).
\]
The balance of sums of the coordinates of weights implies
$n(l+m-t)=p+q-a-2b$ wherefrom either $t=l+m$, $a+2b=p+q$ or
$l=l+m-1$, $a+2b=p+q-n$.

But the second case is impossible because in this case $\nu_n=l+m-3$; hence,
$\lambda_n=l-1$, and therefore $\Pi_if=0$ for $i<n$ (since the last coordinate
of the weight of $\partial_iv'\otimes w'$ cannot be equal to $l+m-3$).
Therefore, $\Pi_{in}f$ is of type 3 or 4) implying
$\lambda_i=\lambda_n=0$,
i.e., $\lambda=(0, \dots , 0, 0)$ and $\lambda_n=l$. This is a~contradiction.

In the first case $t-2=\nu_n=\lambda_n+\mu_n-1=m-2$, hence,
$t=m$, $l=0$. Thus, the first multiple is a~twist-free form.
Such operators will be considered in the next section. Therefore, let us pass to
another case.

Let now $\nu_1-\nu_n\leq 1$, i.e., the image of the operator is also a~twisted form. The fact that
such operators exist only in the cases listed easily follows from the balance
of the sums of coordinates: let
$\nu=(t, \dots , t, \underbrace{t-1, \dots , t-1}_{r})$.

Then, 
$0\leq r\leq n-1$. But
$$
\sum\lambda_i=nl-p, \quad \sum\mu_i=nm-q, \quad \sum\nu_i=nt-r, \quad
\sum\lambda_i+\sum\mu_i-1=\sum\nu_i,
$$i.e.,
$n(l+m-t)=p+q+1-r$, where $-n+2\leq p+q+1-r\leq 2n-1$, and therefore 
either

I) $t=l+m$, $r=p+q+1$ if $p+q+1<n$,
or

II) $t=l+m-1$, $r=p+q+1-n$ if $p+q+1\geq n$.

We have to prove that in each case the operator is unique. To this
end, by the routine method we have to indicate the component which the
singular highest weight vector must contain. If there are two
non-proportional highest weight singular vectors, then certain linear
combination of them does not contain the component indicated which is
impossible.

Let us assume that $l, m\neq 0$ (the opposite case is considered in
the next section). Observe that $\Pi_nf$ must be highest with respect
to $\fgl(E)$, where $E=\langle e_1, \dots , e_{n-1}\rangle$ because the
components $(x_i\partial_j)(\Pi_nf)$, where $i<j<n$ cannot cancel with
the other components of $(x_i\partial_j)f$, and therefore by Lemma \ref{L5.1} the vector $\Pi_nf$ contains the component with vector $v$
highest with respect to $\fgl(E)$. In other words, the weight of is
equal to either $(l, \dots , l, l-1, \dots , l-1)$ or $(l, \dots , l,
l-1, \dots , l-1, l)$.

I) $\nu_n=l+m-1$, hence, $\Pi_nf$ contains summands of the form
$m\partial_nv\otimes w-lv\otimes\partial_nw$, where the $n$-th
coordinate of the weight of $v$ is equal to $l$ that of $w$ is equal
to $m$. Therefore, there should be a~component $lv'_0\otimes\partial
w'$, where the weight of $v'_0$ is equal to $(l, \dots , l, l-1, \dots ,
l-1, l)$. The presence of this component is compulsory.

II) $\nu_nl+m-2$, hence, $\Pi_nf$ may contain summands of the form
$(m-1)\partial_nv\otimes w-lv\otimes\partial_nw$ (hence, the $n$-th
coordinate of the weight of $v$ is equal to $l$ that of $w$ is equal
to $w-m-1$) and of the form $m\partial_nv'\otimes
w'-(l-1)v'\otimes\partial_nw'$ (hence, the $n$-th coordinate of the
weight of $v'$ is equal to $l-1$ that of $w'$ is equal to $m$).

Let us establish when only the summands of the second type may occur. In this
case if $\lambda_i=l$ (such $i$ exist since $p<n$) the vector $\Pi_{in}f$
has no components of the form $p(\partial)v_1\otimes w_0$
(see the list of singular vectors for $n=2$), but for $l\neq0$, $m\neq 0$
this is impossible.

Therefore, this case is excluded and $lv\otimes\partial_nw$ is the
``compulsory'' component; the weight of $v$ is equal to $(l, \dots , l,
l-1, \dots , l-1, l)$. The uniqueness is proved.

\section{The case $T(V_i)=\Omega^p$}

\begin{Lemma}\label{L8} Let $B: (\Omega^p, T(V_1))\to T(V_2)$ be an invariant differential
operator of order $d$. Then, 

{\em 1)} If there exist $\omega\in\Omega^p$ such that $d\omega=0$ and
$B(\omega, s)\not\equiv 0$, then:
when $p=0$, there exists an operator $U: T(V_1)\to T(V_2)$ of order $d$ such that $B(1, s)=U$; and when $p \ge 1$, there exists an invariant
differential operator $B^+: (\Omega^{p-1}, T(V_1))\to T(V_2)$ of order $d+1$
such that $B^+(\omega, s)=B(d\omega, s)$.

{\em 2)} If $d\omega =0$ implies $B(\omega, s)=0$, then there exists $B^-: (\Omega^{p+1}, T(V_1))\to T(V_2)$ an invariant operator of order $d-1$ such that $B(\omega, s)=B^-(d\omega, s)$ (if $d=0$, this case is excluded).
\end{Lemma}

\begin{proof} 1) By definition $B^+(\omega, s)=B(d\omega, s)$ is invariant
since it is expressed in terms of invariant operators. The operator
$B^+$ does not vanish identically since there exists
$\omega_1\in\Omega^p$ such that $d\omega_1=0$ but $B(\omega_1,
s)(x_0)\neq 0$ and there exists $\omega_0\in\Omega^{p-1}$ such that
$d\omega_0$ in a~neighborhood of point $x_0$ such that $B^+(\omega_0,
s)(x_0)=B(\omega_1, s)(x_0)\neq 0$.

2) Let $\omega_1, \omega_2\in\Omega^p$ be such that
$d\omega_1=d\omega_2=\omega\in\Omega^{p+1}$. Then, 
$d(\omega_1-\omega_2)=0$, and therefore $B(\omega_1-\omega_2, s)=0$,
that is, $B(\omega_1, s)=B(\omega_2, s)$.  This shows that: if $\omega=d\omega'$ in a neighborhood of $x$, then $B^-(\omega, s)(x)=B(\omega', s)(x)$ is well-defined; and that, if $d\omega\neq 0$, then  $B^-(\omega, s)=0$.

The
invariance is obvious:
$$
gB^-(\omega, s)=gB(\omega', s)=B(g\omega', gs)=B^-(g\omega, gs).
$$
Let us prove that $B^-$ is local. Let $x_0\not\in\supp\omega$. Then, 
$x_0\in\supp\int^x_{x_0}\omega$,
where $\int^x_{x_0}\omega$ is determined in a~neighborhood of $x_0$ because
$\omega$ is closed.

Hence, $B(\int^x_{x_0}\omega, s)(x_0)=0$, i.e., $B^-(\omega, s)(x_0)=0$.
If $x_0\not\in\supp s$, then
$$
B^-(\omega , s)(x_0)=B(\omega'_, s)(x_0)=0.
$$
The locality is proved. Therefore, $B^-$ is a~differential operator.

But $B(\omega, s)=B^-(d\omega, s)$, and therefore the order of
$B^-$ is equal to $d-1$.
\end{proof}

\begin{Corollary}In $\S 1$, there are listed all the first order operators of the form
$$
B: (\Omega^p, T(V_1))\to T(V_2)
$$
$$
B: (T(V_1), \Omega^q)\to T(V_2)
$$
$$
B: (T(V_1), T(V_2))\to \Omega^r.
$$
\end{Corollary}

\begin{proof} Indeed: let $B:(\Omega^p, T(V_1))\to T(V_2)$ be a~first order
operator. The following cases are possible:

1) There exists a~second order operator $B^+(\Omega^{p-1}, T(V_1))\to
T(V_3)$ such that $B^+(\omega, s)$ is equal to $B(d\omega, s)$. But all such
second order operators are known; these are
$$
\begin{array}{ll}
 B^+_1(\omega_1, \omega_2)&=Z(d\omega_1, d\omega_2), \cr
B^+_2(\omega, s)&=dZ(d\omega, s), \cr
B^+_3(\omega, s)&=P_4(d\omega,
s)\text{ for } p=n-1).
\end{array}
$$
The corresponding first order operators are
$$
\begin{matrix}
B_1(\omega_1, \omega_2)=Z(\omega_1, d\omega_2), \cr
B_2(\omega,s)=dZ(\omega, s),\cr
B_3(\delta, s)=P_4(\delta, s).
\end{matrix}
$$
If $B(\omega, s)-B_j(\omega, s)\neq 0$, then go to case 3).

2) $p=0$ and $B(1, s)U(s)\not\equiv 0$. But all the unary operators
are known: $U(\omega)=d\omega$, $B(1, \omega)=d\omega$, $B(\varphi,
\omega)-\varphi d\omega=0$ if $d\varphi =0$, and therefore $B(\varphi,
\omega)-\varphi d\omega$ which corresponds to case 3).

3) There exists an operator $B^-(\Omega^{p+1}, T(V_1))\to T(V_2)$ of
order zero such that $B(\omega, s)$ is equal to $B^-(d\omega, s)$. But the zero
order operator is $Z(\omega, s)$, hence,
$$
B(\omega, s)=Z(d\omega, s)=P_1(\omega, s).
$$
In the proof of the existence of the operators $B: (T(V_1), \Omega)\to
T(V_2)$ the arguments are the same and as in that of the operators of
the form $B(T(V_1), T(V_2))\to\Omega$, these operators are conjugate
to the already considered operators, and therefore all of them are
listed. Proof is completed.
\end{proof}

\section{Higher order operators}

In this section we will prove that there are no invariant differential
operators of order $\geq 4$.

Let
$$
f=\sum P_i(\partial'{}_1, \dots , \partial''{}_n)U_i, \quad f\in I(V^*_1W^*).
$$
Let $E=\langle E_{i_i}, \dots , e_{i_j}\rangle\subset \Kee^n$ and
$\fgl(E)\subset \fgl(n)$. Recall that
$$
I_E(V^*)=\Kee[\partial'_{i_1}, \dots , \partial'{}_{i_\gamma}]\otimes V^*; \quad
I_E(V^*, W^*)=I_E(V^*)\otimes I_E(W^*).
$$
The vector $f$ can be represented as a~polynomial in
$$
\partial'{}_{j_i},
\dots , \partial'{}_{i_{\bar\gamma}}, \partial''_{j_1}, \dots ,
\partial''_{j_{\bar\gamma}}, \text{ where }\{j_i, \dots ,
j_{\bar\gamma}\}=\{1, \dots , n\}\setminus\{i_j, \dots ,
i_{\gamma}\}
$$
with coefficients from $I_E(V^*)$, namely,
$$
f=\sum P_i(\partial'{}_{j_{1}}, \dots ,
\partial''_{j_{\bar\gamma}}f_{i},\quad f_{i}\in I_E(V^*).
$$
In particular,
the constant term of this polynomial is $\Pi_Ef$.

If $f$ is $\fgl(n)$-highest, then it is also $\fgl(E)$-highest, i.e.,
$x_{\alpha}\partial_{\beta}f=0$ for $\alpha<\beta$ and
$\alpha, \beta\in\{i_1, \dots , i_{\gamma}\}$.
From the equation
$$
x_{\alpha}\partial_{\beta}(\sum P_i(\partial'{}_{j_i}, \dots ,
\partial''_{j_{\bar\gamma}})f_i=\sum P_i(\partial'{}_{j_i}, \dots ,
\partial''_{j_{\gamma}})(x_{\alpha}\partial_{\beta})f_i= 0
$$
it follows that
$(x_{\alpha}\partial_{\beta})f_i=0$, hence, all
the vectors $f_i$ are $\fgl(E)$-highest ones.

Similarly, if $f$ is singular with respect to $\cL$, then all the coefficients
$f_i$ are singular with respect to $\cL_E$.

Now let us pass to the proof (of the fact that there are no invariant
operators of order $>3$).

For $n=1$ the proof was carried out in $\S 3$.

Let now $n=2$. The highest singular vector cannot contain components
in which either $\partial_1$ or $\partial_2$ enters otherwise the
above decomposition contains a~highest weight vector of degree $d\leq
4$ in dimension $n=1$.

Here is the generic form of the vector $f$ of degree 4 without components of the
form $\partial^4_iU$:
$$
f=f_1+f_2+f_3+f_4+f_5,
$$
where
$$
f_1=\partial'{}^3_1\partial'{}_2u_1 + \partial'{}^2_1\partial'{}^2_2u_2 +
 \partial'{}_1\partial'{}^3_2u_3,
$$
$$

\begin{array}{l}
f_2=\partial'{}^3_1\partial''_2u_4 + \partial'{}^2_1\partial'{}_2\partial''_1u_5 +
 \partial'{}^2_1\partial'{}_2\partial''_1u_6 +\\
 \partial'{}_1\partial'{}^2_2\partial''_1u_7 +
 \partial'{}_1\partial'{}^2_2\partial''_2u_8 + \partial'{}^3_2\partial''_1u_9,
\end{array}
$$
$$

\begin{array}{l}
f_3=\partial'{}^2_1\partial''_1\partial''_2u_{10}+
 \partial'{}_1\partial'{}_2\partial''{}^2_1u_{11} +
 \partial'{}^2_1\partial''{}^2_2u_{12} +\\
 \partial'{}_1\partial'{}_2\partial''_1\partial''_2u_{13}+
 \partial'{}^2_2\partial''{}^2_1u_{14} +\\
 \partial'{}_1\partial'{}_2\partial''{}^2_2u_{15} +
 \partial'{}^2_2\partial''_1\partial''_2u_{16}
\end{array}
$$
The form of $f_4$ and $f_5$ is similar to that $f_2$ and $f_1$.

Thus, we have
$$

\begin{array}{l}
X_+f_1=-3\partial'{}^2_1\partial'{}^2_2u_1 - 2\partial'{}_1\partial'{}^3_2u_2 -
 \partial'{}^4_2u_3 + \partial'{}^3_1\partial'{}_2(X_+u_1) +\\
 \partial'{}_1\partial'{}^3_2(X_+u_3);
 \end{array}
$$
$$

\begin{array}{l}
X_+f_2 =-3\partial'{}^2_1\partial'{}_2\partial''_2u_4 -
 2\partial'{}_1\partial'{}^2_2\partial''_1u_5 -
 \partial'{}^2_1\partial'{}_2\partial''_2u_5 -\\
 2\partial'{}_1\partial'{}^2_2\partial''_2u_6 - \partial'{}^3_2\partial''_1u_7 -
 \partial'{}_1\partial'{}^2_2\partial''_2u_7 -\partial'{}^3_2\partial''_2u_8 -
 \partial'{}^3_2\partial''_2u_9 + \\
 \partial'{}^3_1\partial''_2X_+u_4 +
 \partial'{}^2_1\partial'{}_2\partial'{}_1X_+u_5 +
 \partial'{}^2_1\partial'{}_2\partial''_2X_+u_6 +\\
 \partial'{}_1\partial'{}^2_2\partial''_1X_+u_7 +
 \partial'{}_1\partial'{}^2_2\partial''_2X_+u_8 +
 \partial'{}^3_2\partial''_1X_+u_9
\end{array}
$$
$$

\begin{array}{l}
X_+f_3=-2\partial'{}_1\partial'{}_2\partial''_1\partial''_2u_{10} -
 \partial'{}^2_1\partial''{}^2_1u_{10} -\partial'{}^2_2\partial''{}^2_1u_{11} -
 2\partial'{}_1\partial'{}_2\partial''_1\partial''_2u_{11} -\\
 2\partial'{}_1\partial'{}_2\partial''{}^2_2u_{12} -
 \partial'{}^2_2\partial''_1\partial''_2u_{13} -
 \partial'{}_1\partial''_2\partial''{}^2_2u_{13} -\\
 2\partial'{}^2_2\partial''_1\partial''_2u_{14} -
 \partial'{}^2_2\partial''{}^2_2u_{15} -\partial'{}^2_2\partial''{}^2_2u_{16} +
 \partial'{}^2_1\partial''_1\partial''_2X_+u_{10} +\\
 \partial'{}_1\partial'{}_2\partial''{}^2_1X_+u_{11} +
 \partial'{}^2_1\partial''{}^2_2X_+u_{12} +
 \partial'{}_1\partial'{}_2\partial''_1\partial''_2X_+u_{13} +\\
 \partial'{}^2_2\partial''{}^2_1X_+u_{14} +
 \partial'{}_1\partial'{}_2\partial''{}^2_2X_+u_{15} +
 \partial'{}^2_2\partial''_1\partial''_2X_+u_{16}
\end{array}
$$
The form of $X_+f_4$ and $X_+f_5$ is similar to that of $X_+f_2$ and $X_+f_1$.

Since $X_+f=0$, it follows that
$$
\left.\begin{matrix}
u_3=0\cr
2u_2=X_+u_3\cr
3u_1=X_+u_2\cr
0=X_+u_1\cr
 \end{matrix}
\right \} \Longrightarrow u_1=u_2=u_3=0, quad f_1=0.
$$
Similarly, $f_5=0$.

Further, it follows that
$$
\left.\begin{matrix}
u_8+u_9=0\cr
u_7=X_+u_9\cr
2u_6+u_7=X_+u_8\cr
2u_5=X_+u_7\cr
3u_4+u_5=X_+u_6\cr
0=X_+u_5\cr
o=X_+u_4\cr
\end{matrix}\right \}
\Longrightarrow \begin{matrix}
u_8=A,\cr
u_9=-A\cr u_7=-X_+A,\cr
u_6=X_+A,\cr
u_5=-\frac{1}{2}(X_+)^2A,\cr
u_4=\frac{1}{2}(X_+)^2A\cr (X_+)^3A=0.\cr
\end{matrix}
$$
Hence, $f_4=0$.

Finally, it follows that
$$
\left.\begin{matrix}
u_{15}+u_{16}=0\cr
u_{13}+2u_{14}=X_+u_{16}\cr
2u_{12}+u_{13}=X_+u_{15}\cr
u_{11}=X_+u_{14}\cr
2u_{10}+2u_{11}=X_+u_{13}\cr
u_{10}=X_+u_{12}\cr
0=X_+u_{11}\cr
0=X_+u_{10}\cr
\end{matrix}\right \}\Longrightarrow
\begin{matrix}
u_{15}=B, \quad u_{16}=-B, \quad u_{13}=C\cr
u_{14}=\frac{1}{2}(-X_+B-C), \quad
u_{12}=\frac{1}{2}(X_+B-C)\cr
u_{11}=-\frac{1}{2}(X_+)^2B, \quad u_{10}=\frac{1}{2}(X_+)^2B\cr
(X_+)^3B=0, \quad X_+C=0.\cr
\end{matrix}
$$
The singularity condition reads as
$$

\begin{array}{l}
(x^2_2\partial_1)f_2 =-2\partial'{}^3x''_-u_4 -
2\partial'{}^2_1\partial''_1x'{}_-u_5-
2\partial'{}^2_1\partial''_2x'{}_-u_6 - \\
2\partial'{}^2_1\partial'{}_2x''_-u_6 -
4\partial'{}_1\partial'{}_2\partial''_1x'{}_-u_7 +
2\partial'{}^2_2\partial''_1u_7 -\\
4\partial'{}_1\partial'\partial''_2X_-u_8 +
2\partial'{}^2_1\partial''u_8 -
2\partial'{}_1\partial'{}^2_2x''_-u_8 - \\
6\partial'{}^2_2\partial''_1X_-u_9 +
6\partial'{}_1\partial'{}_2\partial''_1u_9
\end{array}
$$
The form of $(x^2_2)f_4$ is similar while
$$

\begin{array}{l}
(x^2_2\partial_1)f_3=-2\partial'{}^2_1\partial''_1x''_-u_{10} -
2\partial'{}_1\partial''{}^2_1x'{}_-u_{11} - 4\partial'{}^2_1\partial''_2x''_-u_{12}+
2\partial'{}_1\partial''_1\partial''_2x'{}_-u_{13} -\\
2\partial'{}_1\partial'{}_2\partial''_1x''_-u_{13} -
4\partial'{}_2\partial''{}^2_1x'{}_-u_{14}+ 2\partial'{}_1\partial''{}^2_1u_{14} -
2\partial'{}_1\partial''{}^2_2X_-u_{15} -\\
4\partial'{}_1\partial'{}_2\partial''_1x''_-u_{15} +
2\partial'{}_1\partial'{}_2\partial''_1u_{15}-
4\partial'{}_2\partial''_1\partial''_2u_{16} +\\
2\partial'{}_1\partial''_1\partial''_2u_{16} - 2\partial'{}^2_2\partial''_1x''_-u_{16}.
\end{array}
$$
Since $(x^2_2\partial_1)f=0$, it follows that
$$

\begin{array}{rl}
x''_-u_4=x''_-u_6=x''_-u_8&=0\\
x'{}_-u_5-u_7+x''_-u_{10}-u_{12}&=0\\
x'{}_-u_6-u_8+2x''_-u_{12}&=0\\
2x'{}_-u_7-3u_9+x'{}_-u_{13}-u_{15}&=0\\
2x'{}_-u_8+2x''_-u_{15}&=0\\
3x'{}_-u_9+x''_-u_{16}&=0
\end{array}
$$
and similar equations that relate $f_3$ with $f_4$.

Since $u_8=-u_9=A$ and $u_{15}=-u_{16}=B$, then the last two equations
imply that
$$
x'{}_-A=x''_-B=0.
$$
But the first equation implies $x''_-A=0$. Moreover, the last two
equations (connecting $f_3$ and $f_4$) imply that
$x'{}_-A=x''_-B=0$. But the first of the equations implies
$x'{}_-B=0$. Thus,
$$
x'{}_-A=x''_-A=0, \quad x'{}_-B=x''_-B=0.
$$
Hence,
$$
A=av_{\lambda}\otimes w_{\mu}, \quad B=bv_{\lambda}\otimes
w_{\mu},
$$
where $v_{\lambda}$ and $w_{\mu}$ are the lowest weight vectors in
$V^*$ and $W^*$, respectively.

But since $(X_+)^3A=0$ and $(X_+)^3B=0$
we deduce that $\lambda +\mu\leq 2$.

Consider the following cases:

1) $f_2\neq 0$ and $f_4\neq 0$. Since $A\neq 0$, the equation
$x''_-u_6=x''_-X_+A$ implies that $\mu =0$.

Similarly, $f_4\neq 0$ implies that $\lambda =0$. But then
$x'{}_-u_6-u_8+2x''_-u_{12}=0$ implies that $u_8=0$, hence, $A=0$.
This is a~contradiction.

2) $f_2=0$. We have
$$

\begin{array}{rlrl}
u_{12}&=x''_-u_{10}&0&=x''_-u_{12}\\
B&=u_{15}=x''_-u_{13}, &B&=x''_-C.\\
B&=(a(02)+b(11)+c(20), &C&=x(01)+y(10).\\
\end{array}
$$

a) $\lambda =0$, $\mu =2$. $a(02)=2x(02)$, $a=2x$. Then, 
$$
u_{12}=\frac{1}{2}(X_+2x(02)-x(01))=\frac{1}{2}x(01), \quad u_{10}=x(00).
$$
The coefficients do not match.

b) $\lambda =\mu =1$, then $B=b(11)=y(11)$ and
$$
\left.\begin{matrix}
u_{15}=-u_{16}=b(11)\cr
u_{12}=b(01), \cr
u_{14}=-b(10), \cr
u_{13}=-b(01)+b(10)\cr
u_{10}=-u_{11}=b(00)
\end{matrix}\right \} \Longrightarrow
\begin{matrix}
u_{12}=\frac{1}{2}(b-x)(01), \quad u_{10}=b(00);\cr
\frac{1}{2}(b-x)=b\Longrightarrow x=-b.
\end{matrix}
$$
We get a~fourth order operator $c(\omega_1,
\omega_2)=d^2\omega_1\wedge d^2\omega_2$ which is invariant only with
respect to $\fsvect(2)$, not $\fvect(2)$.

c) $\lambda =2$, $\mu =0$. Then, $c=y(10)$, $u_{12}=0\Longrightarrow c=0$.
The generic vector of degree 5 is of the form
$$
f=f_1+f_2+f_3+f_4+f_5,
$$
where
$$
f_1=\partial'{}^3_1\partial'{}^2_2u_1 +
\partial'{}^2_1\partial'{}^3_2u_2,
$$
$$

\begin{array}{l}
f_2=\partial'{}^3_1\partial'{}_2\partial''_2u_3 +
 \partial'{}^2_1\partial''_1\partial'{}^2_2u_4 +\\
 \partial'{}^2_1\partial'{}^2_2\partial''_2u_5 +
 \partial'{}_1\partial''_1\partial'{}^3_2u_6
\end{array}
$$
$$

\begin{array}{l}
f_3=\partial'{}^3_1\partial''{}^2_2u_7 +
 \partial'{}^2_1\partial''_1\partial'{}_2\partial''_2u_8+
 \partial'{}_1\partial''{}^2_1\partial'{}_2u_9 +\\
 \partial'{}^2_1\partial'{}_2\partial''{}^2_2u_{10} +
 \partial'{}_1\partial''_1\partial'{}^2_2\partial''_2u_{11} +
 \partial''{}^2_1\partial'{}^3_2u_{12}
\end{array}
$$
and the form of $f_4$, $f_5$ and $f_6$ is similar to that of $f_5$,
$f_2$ and $f_1$, respectively.

We have
$$

\begin{array}{l}
X_+f_1=-3\partial'{}^2_1\partial'{}^2_2u_1 -
2\partial'{}_1\partial'{}^4_2u_2 +\\
 \partial'{}^3\partial'{}^2_2(X_+u_1) +
 \partial'{}^2_1\partial'{}^3_2(X_+u_2);
\end{array}
$$
$$

\begin{array}{l}
X_+f_2 =-3\partial'{}^2_1\partial'{}^2_2\partial''_2u_3 -
 2\partial'{}_1\partial'{}^3_2\partial''_1u_4 -
 \partial'{}^2_1\partial'{}^2_2\partial''_2u_4 -\\
 2\partial'{}_1\partial'{}^3_2\partial''_2u_5 -
\partial'{}^4_2\partial''_1u_6+\text{ components with the } X_+u_i.
\end{array}
$$

Since $X_+f_i=0$, it follows that
$u_1=u_2=0$, $u_6=0$, $2u_5+u_7=0$, $2u_4=X_+u_6, \; (=0)$,
$3u_3+u_4=X_+u_5$, $X_+u_3=X_+u_4=0$, $u_5=A$, $u_7=-2A$,
$u_3=\frac{1}{3}X_+A$,
$(X_+)^2A=0$.
Consider
$$

\begin{array}{l}
X_+f_3=-3\partial'{}^2_1\partial'{}_2\partial''{}^2_2u_7 -
 2\partial'{}_1\partial'{}^2_2\partial''_1\partial''_2u_8 -
 \partial'{}^2_1\partial'{}_2\partial''{}^2_2u_8 -\\
 \partial'{}^3_2\partial''{}^2_1u_9 -
 2\partial'{}_1\partial'{}^2_2\partial''_1\partial''_2u_9 -
 2\partial'{}_1\partial'{}^2_2\partial''{}^2_2u_{10} -\\
 \partial'{}^3_2\partial''_1\partial''_2u_{11} -
 \partial'{}_1\partial'{}^2_2\partial''{}^2_2u_{11} -
 2\partial'{}^3_2\partial''_1\partial''_2u_{12} + \dots ,
 \end{array}
$$
From $X_+f_3=0$ we deduce that
$$

\begin{array}{l}
 2u_{10}+u_{11}=0, \\
 u_{11}+2u_{12}=0, \\
 u_{10}=u_{12}=B, \\
 u_{11}=-2B, \\
 3u_7+u_8=X_+u_{10}, \\
 2u_8+2u_9=X_+u_{11}, \\
 u_9=X_+u_{12}, \\
 u_9=X_+B, \\
 u_8=-2X_+B, \\
 u_7=X_+b, \\
 (X_+)^2B=0, \\
 X_+u_7=X_+u_8=X_+u_9=0.\end{array}
$$
Finally, $f_1=f_2=f_3=0$ and, similarly, $f_4=f_5=f_6=0$.

\subsection*{Acknowledgements} I am thankful to A.~Kirillov for raising the problem,
V.~Palamodov, A.~Kirillov and D.~Leites for general help.

\end{paper}